\documentclass[11pt]{amsart}
\usepackage{amsmath,amssymb,mathrsfs,amscd}
\usepackage{graphicx}
\usepackage[all]{xy}

\title{Double solid twistor spaces: the case of arbitrary signature}

\author{Nobuhiro Honda}
\thanks
{$^{\dag}$This work was partially supported by
Research Fellowships of the 
Japan Society for the Promotion
of Science for Young Scientists.\\
{\it{Mathematics Subject Classifications}} (2000) 32L25, 32G05, 32G07, 53A30,
53C25\\
{\it{Keywords}}\ \  twistor space, minitwistor space, Moishezon manifold, conic bundle, double covering, bimeromorphic transformation,  self-dual metric}
\date{}

\newcommand{\ol}{\overline}
\newcommand{\ra}{\rightarrow}
\newcommand{\lra}{\longrightarrow}

\newcommand{\set}{\,|\,}
\newcommand{\proofend}{\hfill$\square$}

\newtheorem{prop}{Proposition}[section]
\newtheorem{lemma}[prop]{Lemma}

\newtheorem{thm}[prop]{Theorem}
\newtheorem{rmk}[prop]{Remark}

\newtheorem{definition}[prop]{Definition}

\begin{document}

\begin{abstract}
In a recent paper (\cite{Hon07-3}) we constructed a series of new Moishezon twistor spaces  which are a kind of variant of the famous LeBrun twistor spaces.
In this paper we explicitly give  projective models of another series of Moishezon twistor spaces on 
$n\mathbf{CP}^2$ for arbitrary $n\ge 3$, which can be regarded as a
generalization of the twistor spaces of  `double solid type' on $3\mathbf{CP}^2$ studied by Kreu\ss ler, Kurke, Poon and the author.
Similarly to the twistor spaces of `double solid type' on $3\mathbf{CP}^2$, 
 projective models of the present twistor spaces have a natural structure of  double covering of a $\mathbf{CP}^2$-bundle over $\mathbf{CP}^1$.
We explicitly give a defining polynomial of the branch divisor of the double covering, whose restriction to fibers is degree four.
If $n\ge 4$ these are new  twistor spaces, to the best of the author's knowledge.
We also compute the dimension of the moduli space of these twistor spaces.
Differently from \cite{Hon07-3}, the present investigation is based on analysis of pluri-(half-)anticanonical systems of the twistor spaces.
\end{abstract}

\maketitle

\bigskip\noindent
\section{ Introduction}

The Weyl curvature tensor of a Riemannian metric on a manifold   is invariant under changes of the metric by multiplying a function.
When the manifold is more than three-dimensional, it measures how the metric is distant from a conformally flat metric.
If the manifold is four-dimensional and oriented, the Weyl curvature tensor decomposes into two components, called self-dual and anti-self-dual part.
A metric on an oriented four-manifold is called \!{\it self-dual}\, if the anti-self-dual part of the Weyl curvature vanishes identically.
A significant property of self-dual metrics is that, by so called a Penrose correspondence, there naturally associates a very special kind of a complex threefold, called the {\em twistor space}.
As a differential manifold, it is a sphere bundle over the four-manifold.
Although its fibers  are complex submanifolds of the twistor space, their normal bundles are non-trivial but the same as that of line in $\mathbf{CP}^3$.

Nowadays thanks to a theorem of C.\,Taubes \cite{T92},  it is known that there exist a plenty of compact manifolds which admit a self-dual metric.
More precisely, any oriented compact 4-manifold admits a self-dual metric, after gluing  sufficiently many complex projective planes.
Accordingly, there exist a huge number of compact twistor spaces.
N.\,Hitchin \cite{Hi81} showed that compact twistor spaces do not admit a K\"ahler metric, except two well-known examples.
This result looks  made it difficult to investigate compact twistor spaces by usual methods in complex algebraic geometry.
But Y.\,S.\,Poon \cite{P86} and C.\,LeBrun \cite{LB91} changed the scene dramatically by constructing explicit examples of Moishezon twistor spaces.
(A compact complex manifold is called Moishezon if it is bimeromorphic to a projective manifold.)
Later, F.\,Campana and B.\,Kreu\ss ler \cite{CK98} found a new series of Moishezon twistor spaces and described their structure in a quite detailed form.
The base 4-manifolds of all these twistor spaces are $n\mathbf{CP}^2$, the connected sums of $n$ copies of complex projective planes.
This is partly because  by a theorem of Campana \cite{C91}, if a compact four-manifold admits a self-dual metric whose twistor space is Moishezon (or of Fujiki's class $\mathscr C$, more strongly), then the four-manifold is (homeomorphic to) the four-sphere or $n\mathbf{CP}^2$.

In \cite{Hon07-3}, the author explicitly constructed a family of Moishezon twistor spaces with $\mathbf C^*$-action on $n\mathbf{CP}^2$ for arbitrary $n\ge 2$, which seem to be new in the case $n\ge 4$.
The main tool for the construction was the anticanonical system of the twistor spaces.
In fact, the image of its associated meromorphic map is a (singular) rational surface whose structure is independent of the number $n$ of connected sum.
Further, the meromorphic map can be regarded as a quotient map by the $\mathbf C^*$-action.
Namely the image surface was so called a minitwistor space.
By analyzing the structure of the meromorphic quotient map, we constructed the twistor spaces on $n\mathbf{CP}^2$ as explicit bimeromorphic modifications of  (again explicitly constructed) conic bundles over (the minimal resolution of) the minitwistor spaces.

In this paper we provide yet another family of Moishezon twistor spaces on $n\mathbf{CP}^2$ for arbitrary $n\ge 3$.
As those in \cite{Hon07-3} explained above, these twistors have a $\mathbf C^*$-action induced by the $U(1)$-symmetries of the corresponding self-dual metrics.
Further, when $n=3$, they coincide with those in \cite{Hon07-3}.
(Namely they  coincide with non-LeBrun twistor spaces on $3\mathbf{CP}^2$ with $\mathbf C^*$-action studied in \cite{Hon07-2}.)
But when $n\ge 4$, they are new twistor spaces, to the best of the author's knowledge.
Among other properties, the most characteristic feature of the present twistor spaces is that {\em their projective models have a natural structure of double covering over some $\mathbf{CP}^2$-bundle over $\mathbf{CP}^1$, 
as a natural generalization of the double covering structure of the twistor spaces on $3\mathbf{CP}^2$ studied in \cite{Hon07-2, KK92, P92}.}
(On the other hand, the twistor spaces we constructed in \cite{Hon07-3} were also generalization of non-LeBrun twistor spaces on $3\mathbf{CP}^2$ with $\mathbf C^*$-action.
But they are different generalization respecting the structure of minitwistor spaces.)

We outline how these results are obtained, by explaining contents of each sections.
In \S \ref{ss-S} we precisely explain what kind of twistor spaces we shall consider, by specifying the structure of a divisor $S$ which is a member of the half-anticanonical system of the twistor space.
In \S \ref{ss-multanS} we state properties of some pluri-anticanonical systems of the surface $S$.
They will be a basis for the following investigation.
In \S \ref{ss-lst} we consider pluri-half-anticanonical systems on the twistor spaces and provide basic commutative diagrams of meromorphic maps which are indispensable for analyzing the structure of the present twistor spaces.
We also prove the first key result (Lemma \ref{lemma-nt1}) that the sum of some degree-one divisors belongs to a system $|(n-1)F|^{\mathbf C^*}$, the subsystem of $|(n-1)F|$ 
consisting of divisors  defined by $\mathbf C^*$-invariant sections of $(n-1)F$, where $F$ denotes the canonical half-anticanonical bundle of the twistor spaces.
Next in \S \ref{ss-mt} we determine the image of the meromorphic map associated to this system in explicit form (Theorem \ref{thm-mt1}).
In particular, we show that it is a normal rational surface in $\mathbf{CP}^{n+1}$ whose degree is $2(n-1)$.
The meromorphic map can be regarded as a (meromorphic) quotient map of the $\mathbf C^*$-action:
thus the image surface is a minitwistor space.
We also determine its singularities and give their minimal resolutions (Prop.\,\ref{prop-mt2}).
In contrast to the situation we obtained in \cite{Hon07-3},
the structure of these minitwistor spaces depend on $n$, and they form $(n-2)$-dimensional moduli.
When $n=3$, they coincides with the minitwistor space studied in \cite{Hon-p-1} and \cite{Hon07-3}.
When $n\ge 4$, they are new minitwistor spaces, as far as the author knows.

In Section 3 by analyzing the meromorphic quotient map in the previous section we realize  projective models of the twistor spaces as  conic bundles over  (the minimal resolution of) the minitwistor spaces.
In \S \ref{ss-pr} we concretely give a partial elimination of the indeterminacy locus of the meromorphic map.
Then in \S \ref{ss-dl} and \S \ref{ss-nt} we explicitly construct a $\mathbf{CP}^2$-bundle over the resolution of the  minitwistor space and realize a projective model of the twistor space as a conic subbundle of this $\mathbf{CP}^2$-bundle, by explicitly giving its defining equation.
The idea used here is similar to those in \cite{Hon07-3}, Section 3.
An interesting observation is that, as in \cite{Hon07-3}, a discriminant curve of the conic bundle is a hyperplane section of the minitwistor spaces with respect to its natural realization in a projective space (Lemma \ref{lemma-hps1}).
We show that the inverse image of this discriminant curve splits into two irreducible components and correspondingly we obtain a pair of (mutually conjugate) $\mathbf C^*$-invariant divisors $Y$ and $\ol{Y}$ in the twistor space (Prop.\,\ref{prop-Y1}).
We next determine the cohomology classes of $Y$ and $\ol{Y}$ (in $H^2(Z,\mathbf Z)$) by showing that  adding some divisor to $Y$ (and $\ol{Y}$) gives  members of the system $|(n-2)F|$ (Prop.\,\ref{prop-nt4}).
Existence of these members is a key in obtaining the following presentations of the twistor spaces.

In Section 4 we investigate a complete linear system $|(n-1)F|$.
Based on the results of the previous section,  especially the existence of the above divisors in  the system $|(n-2)F|$, we concretely give generators of $|(n-1)F|$.
Then using them we explicitly determine the image of the meromorphic map associated to $|(n-1)F|$ (Theorem \ref{thm-bim}).
We also show that the map is bimeromorphic onto the image.
In particular, we obtain the second projective models of the twistor spaces.
They are birational to codimension two subvarieties in certain $\mathbf{CP}^4$-bundle over $\mathbf{CP}^1$, whose restriction to each fibers are quartic surfaces.

In the final section, combining results obtained so far, we investigate the meromorphic map associated to the system $|(n-2)F|$.
We show that the map is generically 2 to 1 over its image, and  the image is a rational scroll of planes in $\mathbf{CP}^n$, whose degree is $n-2$.
By blowing-up its vertices which form a line, we obtain a $\mathbf{CP}^2$-bundle $\mathbf P(\mathscr O(n-2)^{\oplus 2}\oplus \mathscr O)$ over $\mathbf{CP}^1$.
Then we explicitly give a meromorphic map from the projective model in Section 4 to this $\mathbf{CP}^2$-bundle, and show that the map is generically 2 to 1.
We further  determine the defining equation of the branch locus of this map in an explicit form.
These results are summarized as Theorem \ref{thm-dc2}.
When $n=3$, the equation of the branch divisor coincides with the one we obtained in \cite{Hon07-2} for a non-LeBrun twistor space on $3\mathbf{CP}^2$ with $\mathbf C^*$-action.
Thus it becomes apparent that {\em the present twistor spaces are natural generalization of the twistor spaces on $3\mathbf{CP}^2$ studied in \cite{Hon07-2}, respecting the structure of double covering.}
We also show that the branch divisor of the double covering is irreducible, non-normal surface, and birational to ruled surface of genus $[(n-1)/2]$ (Prop.\,\ref{prop-B}).
We further obtain some constraint on the defining equation of the branch divisor (Prop.\,\ref{prop-dr}), which is a generalization of the constraint appeared in the case of $n=3$ in \cite{Hon07-2}.
Consequently  the moduli space of the present twistor spaces can be computed to be $n$-dimensional.
Finally we remark that the present twistor spaces can be obtained as a small deformation of LeBrun metrics with torus action which preserves a particular $U(1)$-action, and compare with the results we obtained in \cite{Hon07-1}.

In summary, we investigate the twistor spaces by using three linear systems:
(a) the non-complete system $|(n-1)F|^{\mathbf C^*}$ from which we derive a conic bundle description over the minitwistor spaces (these are studied in Sections 2 and 3),
(b) the complete linear system $|(n-1)F|$ which gives projective models that are  birational to the total space of a fiber space over  $\mathbf{CP}^1$ whose fibers are quartic surfaces
(these are studied in Section 4),
(c) the complete linear system $|(n-2)F|$ which gives a generically 2 to 1 covering over a rational scroll of planes in $\mathbf{CP}^n$ (these are studied in Section 5).
Of course, the last presentation is most simple because the double covering is uniquely determined by only specifying the branch divisor which is expressed by a {\em single }polynomial.
It seems difficult for the author to obtain the last presentation without studying both (a) and (b) in detail (as we did in this paper).

\bigskip\noindent
{\bf Notations and Conventions.}
As in \cite{Hon07-3}, to save  notations we adapt the following convention.
If $\mu:X\to Y$ is a bimeromorphic morphism of  complex variety and $W$ is a complex subspace in $X$, we write $W$ for the image $\mu(W)$  if the restriction $\mu|_W$ is still bimeromorphic.
Similarly, if $V$ is a complex subspace of $Y$, we write $Y$ for the strict transform of $Y$.
If $D$ is a divisor on a variety $X$, the dimension of a complete linear system $|D|$ always means $\dim H^0(X,[D])-1$.
The base locus is denoted by ${\rm Bs}\,|D|$.
If a Lie group $G$ acts on $X$ by means of biholomorphism and $D$ is $G$-invariant,
$G$ naturally acts on the vector space $H^0(X,[D])$.
Then $H^0(X,[D])^{G}$ means the subspace of $G$-invariant sections.
Further $|D|^G$ means its associated linear system.
If $s$ is a non-zero section of a holomorphic line bundle on $X$, $(s)$ denotes the divisor defined as the zero locus of $s$.
When discussing cohomology classes represented by  complex curves $C_1$ and $C_2$ on a complex surface $S$,
we write $C_1\sim C_2$ to mean that they are cohomologous; namely if  $C_1$ and $C_2$ determine the same element in $H^2(S,\mathbf Z)$.
(This is used only in Section 3.)

If $Z$ is a twistor space, $F$ always denotes the canonical square root of the anticanonical line bundle of $Z$ (often called the `fundamental line bundle').
The degree of a divisor on $Z$ means its intersection number with twistor lines.

\section{Analysis of the structure of minitwistor spaces}
\subsection{Construction of the surface $S$ contained in the fundamental system}
\label{ss-S}
As in \cite{Hon07-3}, we start our investigation by specifying structure of a real irreducible member of the fundamental system on  the twistor spaces.

So we consider the product surface $\mathbf{CP}^1\times\mathbf{CP}^1$ equipped with a real structure given by 
\begin{align}\label{rs1}
\text{(complex conjugation)}\times\text{(anti-podal)}.
\end{align}
We write $\mathscr O(1,0)=p_1^*\mathscr O(1)$ and $\mathscr O(0,1)=p_2^* \mathscr O (1)$, where $p_i$ denotes the projection to the $i$-th factor.
We choose a non-real curve $C_1\in|\mathscr O(1,0)|$ and a point $P_1\in C_1$.
Next we blow-up $\mathbf{CP}^1\times\mathbf{CP}^1$ at $P_1$ and $\ol{P}_1$ to obtain a surface with $c_1^2=6$, which has distinguished curves $C_1$ and $\ol{C}_1$ satisfying $C_1^2=\ol{C}_1^2=-1.$
Subsequently we blow-up the surface at the intersection points of $C_1\cup\ol{C}_1$ and the exceptional curves of the last blow-ups.
Then we obtain a surface with $c_1^2=4$ possessing distinguished curves $C_1$ and $\ol{C}_1$ satisfying $C_1^2=\ol{C}_1^2=-2$.
By repeating this procedure $(n-1)$-times, we obtain a surface $S'$ with $c_1^2=8-2(n-1)$ possessing distinguished curves $C_1$ and $\ol{C}_1$ satisfying $C_1^2=\ol{C}_1^2=1-n$.
This surface $S'$ has an obvious structure of  toric surface.
Then as a final operation to obtain the required surface, we choose a conjugate pair of points which are on the exceptional curves of the final blow-ups (in obtaining the surface $S'$) but which are {\em not}  fixed points of the torus action.
Then the action of the torus is killed and 
consequently we obtain a non-toric surface $S$ with a non-trivial $\mathbf  C^*$-action
which satisfies $c_1^2=8-2n$.
The surface $S$ has a natural real structure induced from \eqref{rs1}.
Let $B_1$ and $\ol{B}_1$ be the exceptional curves of the final blowing-up $S\to S'$.
By choosing the pair of blown-up point generically,
we can suppose that there exists no $\mathbf C^*$-invariant $(-2)$-curve intersecting both of $B_1$ and $\ol{B}_1$.
Namely we suppose that the pair of blown-up points on $S'$ do not belong to the same $\mathbf C^*$-orbit closure.
Then there exists a unique pair of $(-1)$-curves $B_2$ and $\ol B_2$ intersecting $\ol B_1$ and $B_1$ transversally respectively.
We note that only freedom involved in the construction of the surface $S$ is the choice of the final blown-up points on $S'$.

\begin{figure}
\includegraphics{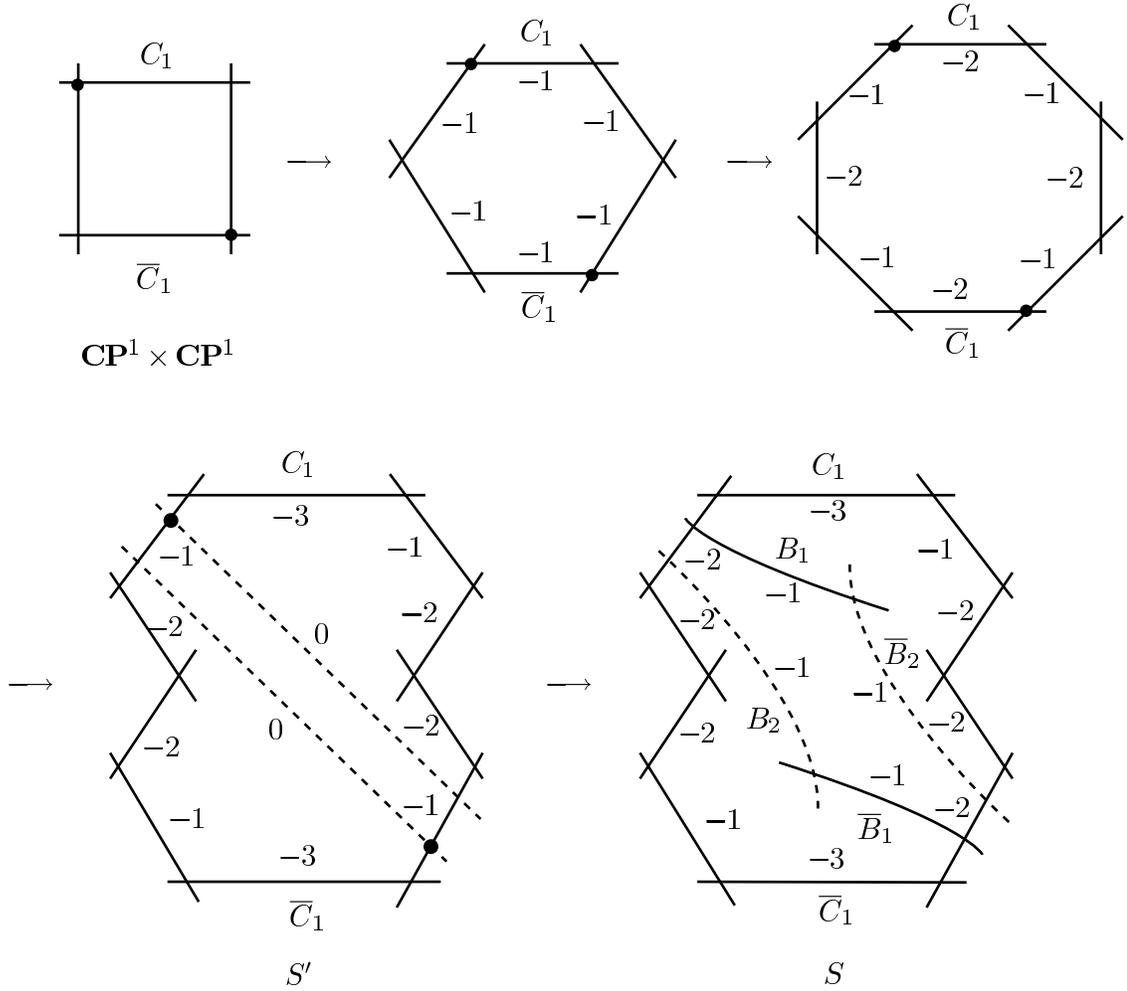}
\caption{Construction of the surface $S$ in the case $n=4$.}
\label{fig-cycle1}
\end{figure}

\begin{figure}
\includegraphics{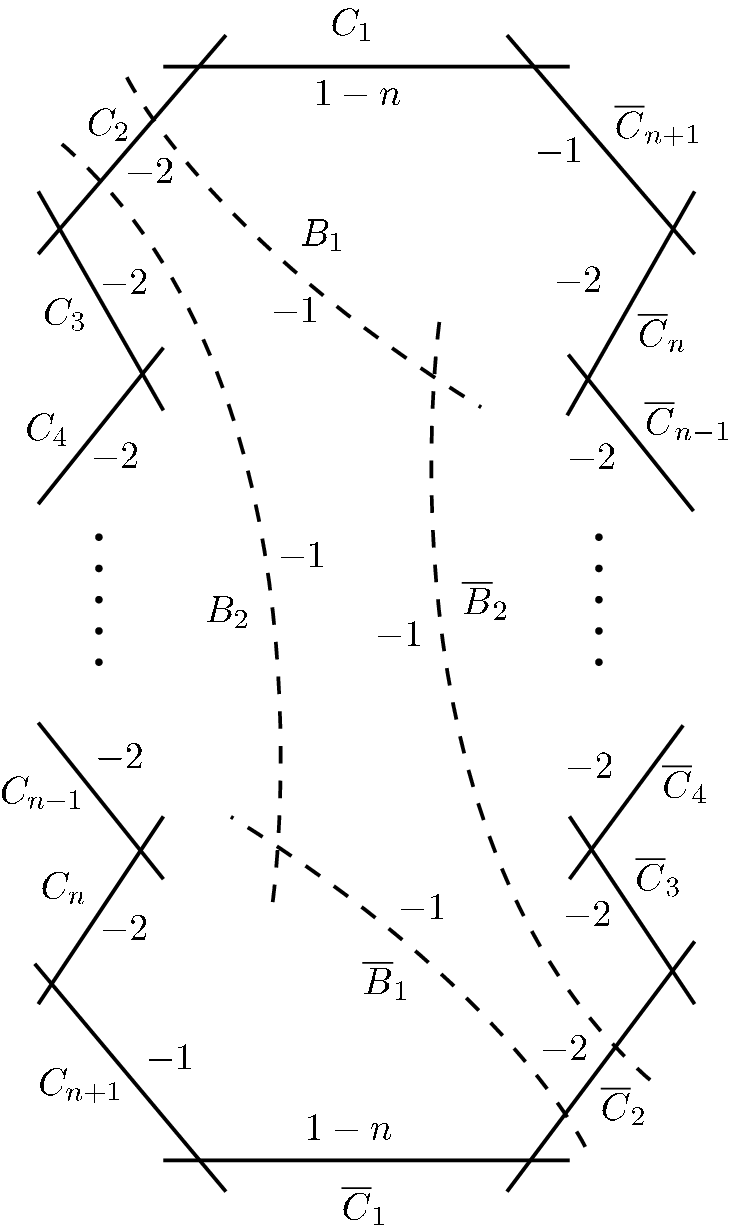}
\caption{The anticanonical cycle in the surface $S$.}
\label{fig-cycle2}
\end{figure}

It is readily seen that if $n\ge 4$ the anticanonical system of the surface $S$ consists of a unique member and it is a cycle of $(2n+2)$ smooth rational curves containing all exceptional curves in obtaining $S'$.
We  write $C$ for this anticanonical cycle (analogously to the notation in \cite{Hon07-3}) and let
\begin{align}
C=\sum_{i=1}^{n+1} C_i+\sum_{i=1}^{n+1}\ol{C}_i
\end{align}
be its decomposition into  irreducible components orderly named as in Figure \ref{fig-cycle2}.
Here, among these, only $C_2$ and $\ol{C}_2$ are  $\mathbf C^*$-fixed.
Note that the pair $B_1$ and $\ol{B}_1$ of  the exceptional curves of the final blow-ups $S\to S'$ are not contained in this cycle $C$.
These intersect $C_2$ and $\ol C_2$ respectively.
By construction
it is easy  to see that the self-intersection numbers of the irreducible components are given by
\begin{align}
C_1^2=1-n,\,\,\, C_2^2=C_3^2=\cdots=C_n^2=-2,\,\,\, C_{n+1}^2=-1.
\end{align}
Of course, we have $\ol{C}_i^2=C_i^2$ for any $i$.
Note that if $n=3$ the construction of the  surface $S$ is the same for that of the surface $S$ treated in \cite{Hon07-3}.
If $n\ge 4$ the present surface is different from \cite{Hon07-3}, since there the unique anticanonical curve was the cycle consisting of 8 irreducible components only.

\subsection{Properties of some pluri-anticanonical systems of $S$}\label{ss-multanS}
We wish to study  algebro-geometric structure of twistor spaces which contain the surface $S$ as its real irreducible member of the fundamental system $|F|$.
Before doing this, in this subsection, we collect basic properties of some pluri-anticanonical systems of $S$.
They will be a basis of our analysis on the multi-systems $|mF|$ on the twistor spaces.

\begin{prop}\label{prop-multan1} Suppose $n\ge 4$ and 
let $S$ be the surface constructed above, and $C$ the anticanonical cycle whose irreducible components are named orderly as in Figure \ref{fig-cycle2}. Then we have the following.
(i) $\dim H^0(mK^{-1}_S)=1$ for any $1\le m\le n-3$.
(ii) $\dim H^0((n-2)K^{-1}_S)=3$ and $\dim H^0((n-2)K^{-1}_S)^{\mathbf C^*}=1$.
(iii) The fixed component of the net $|(n-2)K^{-1}_S|$ is given by
\begin{align}\label{fixcompo1}
(n-3)\sum_{i=1}^4(C_i+\ol{C}_i)+\sum_{i=5}^n (n+1-i)(C_i+\ol{C}_i).
\end{align}
(iv) The  movable part of the net $|(n-2)K_S^{-1}|$ does not have base points.
(v) The associated morphism $S\to\mathbf{CP}^2$ is the composition of the following 2 morphisms: the first one is a succession of blowing-downs of $C_{n+1}, C_n,\cdots, C_6, C_5$ and their conjugates, and the second one is a generically 2 to 1 map from the resulting surface to $\mathbf{CP}^2$. 
Further, the branch curve of the second map is a sum of two different $\mathbf C^*$-invariant irreducible conics, where $\mathbf C^*$ acts on $\mathbf{CP}^2$ by 
\begin{align}\label{action49}
(y_0,y_1,y_2)\longmapsto (y_0,sy_1,s^{-1}y_2), \,\,\,s\in\mathbf C^*
\end{align}
in a homogeneous coordinate.
\end{prop}

Although these are not straightforward to see, 
we omit a proof since they can be shown by standard arguments and computations in algebraic geometry.
Note that by (v) the morphism $S\to\mathbf{CP}^2$ associated to the movable part of the system $|(n-2)K_S^{-1}|$ is generically 2 to 1, and its discriminant locus, consisting of two $\mathbf C^*$-invariant conics,  has two $A_3$-singularities at the intersection points.
Hence the double covering in the usual sense has two $A_3$-singularities over there.
The components $C_1,C_2, C_3$ and their conjugates arise as the exceptional curves of their minimal resolutions, and 
$C_4$ and $\ol{C}_4$ are mapped biholomorphically to the line $\{y_0=0\}$ in the coordinate of \eqref{action49}.
Thus the morphism $S\to\mathbf{CP}^2$ factors as
\begin{equation}\label{factor1}
S\lra \tilde{S}_0\lra S_0\lra \mathbf{CP}^2,
\end{equation}
where $S_0\to\mathbf{CP}^2$ is the double covering branched along
the sum of two $\mathbf C^*$-invariant conics, $\tilde{S}_0\to S_0$ is the minimal resolution of the two $A_3$-singularities, and 
$S\to \tilde{S}_0$ is the contraction of $C_{n+1}, C_n,\cdots, C_6, C_5$ and their conjugates.
Thus the structure of $\tilde S_0$ is independent of $n$.
Note that by (ii) of Prop.\,\ref{prop-multan1} non-zero section of $H^0((n-2)K_S^{-1})^{\mathbf C^*}$ is unique up to constant.
Its zero divisor is exactly the non-reduced curve $(n-2)C$.

Note also that Prop.\,\ref{prop-multan1} holds even when $n=3$.
In that case (i)  claims nothing and (ii) means $|K_S^{-1}|$ is a net.
Further
 every components of \eqref{fixcompo1} disappear (since $C_5,C_6,\cdots$ do not exist) and hence the anticanonical system becomes free.
Moreover $S=\tilde{S}_0$ holds in \eqref{factor1}.
This is exactly the situation we encountered in \cite{Hon07-2}, where we used the double covering map to detect  twistor lines.

We will also need the following results concerning the system $|(n-1)K^{-1}_S|$.

\begin{prop}\label{prop-multan2}
Let $S$ be as in Prop.\,\ref{prop-multan1}. Then we have the following.
(i) $\dim H^0((n-1)K^{-1}_S)=5$ and $\dim H^0((n-1)K^{-1}_S)^{\mathbf C^*}=3$.
(ii) The fixed component of the system $|(n-1)K^{-1}_S|$ is given by
\begin{align}
(n-2)\left(C_1+\ol{C}_1+C_2+\ol{C}_2\right)+\sum_{i=3}^n(n-i+1)(C_i+\ol{C}_i).
\end{align}
(iii) The movable part of the 4-dimensional system $|(n-1)K^{-1}_S|$ does not have base point.
(iv) The associated morphism $S\to\mathbf{CP}^4$ is birational onto its image, and the image is a quartic surface.
(v) The meromorphic map associated to the net $|(n-1)K_S^{-1}|^{\mathbf C^*}$ is a morphism and its image is a conic in $\mathbf{CP}^2$.
Further, its general fibers are smooth rational curves intersecting $C_2$ and $\ol{C}_2$.
\end{prop}

Again we do not write a proof of this proposition.
In short, these two propositions mean that the system $|(n-2)K_S^{-1}|$ gives a generically 2 to 1 morphism onto $\mathbf{CP}^2$ and $|(n-1)K^{-1}_S|$ gives a birational morphism onto a quartic surface in $\mathbf{CP}^4$.


\subsection{Analysis of linear systems on the twistor spaces}\label{ss-lst}
With these preliminary results on the structure of the surface $S$, we begin our study on twistor spaces. 
The following proposition is about the structure of the system $|F|$ and easy to see.

\begin{prop}\label{prop-fundamental}
Suppose $n\ge 4$ and let $Z$ be a twistor space on $n\mathbf{CP}^2$ with $\mathbf C^*$-action which is compatible with the real structure.
Suppose that the complex surface $S$ in Prop.\,\ref{prop-multan1}( and \ref{prop-multan2})
is equivariantly contained in $Z$ as a real $\mathbf C^*$-invariant member of the system $|F|$.
Then we have the following.
(i) $\dim |F|=\dim |F|^{\mathbf C^*}=1$ and {\rm Bs}\,$|F|=C$.
(ii) The pencil $|F|$ has precisely $(n+1)$ reducible members, and they are of the form $S_i^++S_i^-$ $(1\le i\le n+1$), where $S_i^+$ and $S_i^-$ are mutually conjugate, $\mathbf C^*$-invariant divisors  of  degree one.
(iii) $L_i:=S_i^+\cap S_i^-$ $(1\le i\le n+1)$ are $\mathbf C^*$-invariant twistor lines joining conjugate pairs of singular points of the cycle $C$.
\end{prop}

We name the reducible members  $S_i^++S_i^-$ $(1\le i\le n+1$) by imposing that the invariant twistor line $L_i=S_i^+\cap S_i^-$ intersects $C_{i-1}\cap C_i$, where $C_0=\ol{C}_{n+1}$.
Further, we distinguish the 2 irreducible components of $S_i^++S_i^-$ by declaring that $S_i^+$ contains the component $C_{n+1}$.

The existence of the surface $S$  in the twistor space has the following consequence on the induced $U(1)$-action on the base 4-manifold.

\begin{prop}\label{prop-u1}
Let $Z$ and $S$ be as in Prop.\,\ref{prop-fundamental}.
Then the induced $U(1)$-action on $n\mathbf{CP}^2$ satisfies the following.
(i) $U(1)$-fixed locus consists of one sphere and $n$ isolated points.
(ii) There exists a unique isolated $U(1)$-fixed point whose twistor line is fixed by $\mathbf C^*$-action.
(iii) There exists a unique isolated $U(1)$-fixed point whose isotropy subgroup of its twistor line is a cyclic subgroup of $U(1)$ (or $\mathbf C^*$) with order $(n-2)$.
(iv) The isotropy subgroup of all other invariant twistor lines are identity.
\end{prop}

\noindent
Proof.
The $U(1)$-fixed sphere in (i) is the image of the $\mathbf C^*$-fixed curve $C_2$ by the twistor fibration $Z\to n\mathbf{CP}^2$.
The twistor line in (ii) is the one going through a $\mathbf C^*$-fixed point $B_1\cap B_2$, where $B_1$ and $B_2$ are $(-1)$-curves in $S$ obtained in the construction of $S$.
The twistor line in (iii) is $L_1$.
Since the intersection point $C_i\cap C_{i+1}$ ($3\le i\le n$ ) are isolated $\mathbf C^*$-fixed points on $S$, $L_i$ ($4\le i\le n+1)$ are over the isolated $U(1)$-fixed points on $n\mathbf{CP}^2$.
Thus we obtain $n$ isolated $U(1)$-fixed points as in (i).
(iv) can also be verified by computing the induced $U(1)$-action on the tangent spaces of $n\mathbf{CP}^2$ at the $U(1)$-fixed points and then computing the induced action on the spin bundles.
\proofend

\bigskip
The twistor line in (ii) will play an important role in our analysis of the structure of the twistor spaces. So we introduce the following.

\begin{definition} \label{def-ftl}
{\em 
We will call the unique twistor line in (ii) {\em the fixed twistor line}.
We always denote it by $L_0$.
}
\end{definition}

Since $mF|_S\simeq mK_S^{-1}$ by adjunction formula, 
Prop.\,\ref{prop-multan1} (i) means that 
in order to obtain (bimeromorphic)  projective models of the twistor spaces,
we need to consider a linear system $|mF|$ for $m\ge n-2$.
To study these systems systematically, we introduce the following.

\begin{definition} \label{def-trivial}(i) 
For each positive $m\in\mathbf Z$, let $V_m\subset H^0(mF)$ be a linear subspace generated by the image of  a natural multi-linear map
\begin{equation}
H^0(F)\times H^0(F)\times\cdots\times H^0(F)\,\lra\, H^0(mF)
\end{equation}
given by $(s_1,s_2,\cdots,s_m)\mapsto s_1\otimes s_2\otimes \cdots\otimes s_m$.
Prop.\,\ref{prop-multan1} (i) and Prop.\,\ref{prop-fundamental} (i) implies $\dim V_m=m+1$ and $V_m\subset H^0(mF)^{\mathbf C^*}$.
(ii) A divisor $D\in|mF|$ is said to be a trivial member if $D\in |V_m|$.
In other words, $D$ is said to be a trivial member if $D=S_1+\cdots+S_m$ for some $S_i\in |F|$. ($S_i$ is not necessarily irreducible, of course.)
(iii) The meromorphic map from $Z$ to $\mathbf P^{\vee}V_m= \mathbf{CP}^m$ associated to the system $|V_m|$ is denoted by $\Psi_m$. 
The image of $\Psi_m$ is denoted by $\Lambda_m$.
$\Lambda_m$ is a rational normal curve. In particular, its degree in $\mathbf{CP}^m$ is $m$.
\end{definition}

Since we consider the systems $|mF|$ for different values of $m$, we introduce the following.

\begin{definition}\label{def-asso1}
For each positive $m\in\mathbf Z$, let 
\begin{align}
\Phi_m:Z\to\mathbf P^{\vee}H^0(mF)\,\,\,{\text{and}}\,\,\,\,\Phi_m^{\mathbf C^*}:Z\to\mathbf P^{\vee}H^0(mF)^{\mathbf C^*}
\end{align}
 be the meromorphic maps associated to the system $|mF|$ and $|mF|^{\mathbf C^*}$ respectively.
\end{definition}

We note an obvious relation between two meromorphic maps $\Psi_m$ and $\Phi_m$ in the two definitions.
Namely for each $m$ there is a commutative diagram of meromorphic maps
\begin{equation}\label{cd1}
 \CD
Z@>{\Phi_m}>>\mathbf P^{\vee}H^0(mF) \\
 @V\Psi_m VV @VV{\pi_m}V\\
\Lambda_m@>>>\mathbf P^{\vee}V_m\\
 \endCD
 \end{equation}
 where $\pi_m$ is the projection induced from the inclusion $V_m\subset H^0(mF)$.
 (The bottom row is an inclusion as a rational normal curve as in Def.\,\ref{def-trivial}.)
We note that the diagram \eqref{cd1} is $\mathbf C^*$-equivariant, where $\mathbf C^*$-actions are trivial on $\Lambda_m$ and $\mathbf P^{\vee}V_m$.
Taking the $\mathbf C^*$-fixed part, we obtain the diagram
\begin{equation}\label{cd2}
 \CD
Z@>{\Phi_m^{\mathbf C^*}}>>\mathbf P^{\vee}H^0(mF)^{\mathbf C^*} \\
 @V\Psi_m VV @VV{\pi_m^{\mathbf C^*}}V\\
\Lambda_m@>>>\mathbf P^{\vee}V_m\\
 \endCD
 \end{equation}
We use these diagrams to reveal a structure of the twistor spaces.
As an immediate consequence of these diagrams, we obtain that the images  $\Phi_m(Z)$ and $\Phi_m^{\mathbf C^*}(Z)$ are always contained in $\pi_m^{-1}(\Lambda_m)$ and 
 $(\pi_m^{\mathbf C^*})^{-1}(\Lambda_m)$ respectively.

Prop.\,\ref{prop-fundamental} and Prop.\,\ref{prop-multan1} (i) imply that if $m<n-2$ we have the coincidence $V_m=H^0(mF)=H^0(mF)^{\mathbf C^*}$.
Hence the projection $\pi_m$ in the diagram \eqref{cd1} is the identity map and therefore the image $\Phi_m(Z)$ is the rational normal curve $\Lambda_m$.
Hence the map $\Phi_m$ gives little information  for $m<n-2$.
If $m=n-2$, Prop.\,\ref{prop-multan1} and the diagram \eqref{cd1} give a possibility that $\Phi_{n-2}$  gives a generically 2 to 1 covering onto its image.
Also, if $m=n-1$, Prop\,\ref{prop-multan2} and the diagram \eqref{cd1} give a possibility that  the map $\Phi_{n-1}$ is bimeromorphic over its image.
(Both of these expectations will turned out to be true.)

A basic tool relevant to this kind of problem is the short exact sequence
\begin{align}\label{ses1}
0\lra (m-1)F\lra mF\lra mK_S^{-1}\lra 0.
\end{align}
At least for the case  $n=4$, this sequence is enough for proving that $\Phi_{n-2}$ is generically 2 to 1 onto its image.
Actually, when  $n=4$, we put $m=n-2=2$ in  \eqref{ses1}.
 In this case, we use  Riemann-Roch formula and Hitchin's vanishing theorem to deduce
 $\dim H^0(F)-\dim H^1(F)=10-2n=2$.
 It then follows  $H^1(F)=0$ by Prop.\,\ref{prop-fundamental} (i).
Hence we obtain the surjectivity of the restriction map $H^0(2F)\to H^0(2K_S^{-1})$.
Therefore  by Prop.\,\ref{prop-multan1}, $S\in|F|$ is mapped surjectively to $\mathbf{CP}^2$  by $\Phi_{n-2}=\Phi_2$ and it is generically 2 to 1.
Hence by the diagram \eqref{cd1}, the image $\Phi_2(Z)$ is 3-dimensional and $\Phi_2$ is generically 2 to 1 onto its image.

When $n\ge 5$ we still put $m=n-2$ in \eqref{ses1}  and it becomes
\begin{equation}\label{ses3}
0\lra (n-3)F\lra (n-2)F\lra (n-2)K^{-1}_S\lra 0.
\end{equation}
However, a computation similar to the above shows that $H^1((n-3)F)\neq 0$ (or more precisely its dimension is  quite  high) and the cohomology exact sequence of \eqref{ses3} does not mean the existence of a non-trivial member of $|(n-2)F|$. 
Namely, if $n\ge 5$ the argument using cohomology exact sequence fails.
This difficulty  always happens when we consider twistor spaces on $n\mathbf{CP}^2$, $n\ge 5$.
For example, in \cite{Hon07-3}, we needed to show that the system $|2F|=|K^{-1}_Z|$ has a non-trivial member, in order to show that the image of  its associated map is 2-dimensional.
This was shown by proving that a sum of some four degree-one divisors gives a non-trivial member.
Unfortunately, in the present case, a computation shows that if  a sum of degree-one divisors belongs to $|(n-2)F|$, then it must be a trivial member, already for the case $n=4$.
So the method of \cite{Hon07-3} does not work directly.

We overcome this difficulty by considering a system of higher degree, $|(n-1)F|$.
For this system we can find a non-trivial member which is a sum of degree-one divisors as follows.

\begin{lemma}\label{lemma-nt1}
Let $Z$ be as in Prop.\,\ref{prop-fundamental}.
Then the $\mathbf C^*$-invariant divisors
\begin{equation}\label{nt1}
(n-2)S_1^+\,+\sum_{i=2}^{n+1} S_i^+\,\,\,\,{\text{and}}\,\,\,\,\,
(n-2)S_1^-\,+\sum_{i=2}^{n+1} S_i^-
\end{equation}
are non-trivial members of the system $|(n-1)F|^{\mathbf C^*}$.
\end{lemma}

\begin{rmk}
{\em
In general, even if a sum of degree-one divisors belongs to a system $|mF|$ for some $m\in\mathbf Z$, it does necessarily belong to the subsystem $|mF|^{\mathbf C^*}$.
}
\end{rmk}

\noindent Proof of Lemma \ref{lemma-nt1}.
To show that the divisors \eqref{nt1} are members of $|(n-1)F|$, 
since the restriction map of the cohomology groups
$H^2(Z,\mathbf Z)\to H^2(S,\mathbf Z)$ is always injective,
it suffices to show that the restriction of the divisors \eqref{nt1} to $S$ belongs to 
$|(n-1)K_S^{-1}|$.
Since the restrictions $S_i^+|_S$ and $S_i^-|_S$ are precisely two halves of the anticanonical cycle $C$ divided by the twistor line $L_i$, 
we can write these restrictions as  sums of irreducible components of $C$.
Once this is done, it is a routine work to show that they belong to $|(n-1)K_S^{-1}|$.
(For example, it suffices to verify that the intersection numbers with generators of $H^2(S,\mathbf Z)$ coincide.)
We omit the detail of the verification.
It can also be verified that the restrictions moreover belong to the subsystem $|(n-1)K_S^{-1}|^{\mathbf C^*}$.
This means that the 2 divisors \eqref{nt1} are contained in the subsystem $|(n-1)F|^{\mathbf C^*}$.

Finally we show that the divisors are non-trivial members.
 If the former divisor of \eqref{nt1} is a trivial member, there must exist $1\le i,j\le n+1$ such that  $S_i^++S_j^+$ is a member of the pencil $|F|$.
 This contradicts Prop.\,\ref{prop-fundamental} (ii).
 Hence the divisor is not a trivial member.
Then by the reality of the system $|V_m$,  the latter divisor is not a trivial  member as well.
  \proofend

\bigskip 
By using Lemma \ref{lemma-nt1} we can explicitly give generators of the system $|(n-1)F|^{\mathbf C^*}$ as follows.

\begin{prop}\label{prop-generator1}
Let $Z$ be as in Prop.\,\ref{prop-fundamental}.
Then we have the following.
(i) $\dim H^0((n-1)F)^{\mathbf C^*}=n+2$.
(ii) As generators of the system $|(n-1)F|^{\mathbf C^*}$ we can take the following $\mathbf C^*$-invariant divisors:
\begin{align}\label{trivial1}
(n-1-k)(S_1^++S_1^-)+k(S_2^++S_2^-),\,\,0\le k\le n-1,
\end{align}
\begin{align}\label{nt2}
(n-2)S_1^+\,+\sum_{i=2}^{n+1} S_i^+,\,\,
(n-2)S_1^-\,+\sum_{i=2}^{n+1} S_i^-.
\end{align}
(iii) ${\rm Bs}\,|(n-2)F|^{\mathbf C^*}=C-C_{n+1}-\ol{C}_{n+1}.$
\end{prop}

Note that the $n$ divisors \eqref{trivial1} are (independent) generators of the $(n-1)$-dimensional linear system $|V_{n-1}|$.

\bigskip
\noindent
Proof of Prop.\,\ref{prop-generator1}.
We first show that $\dim H^0((n-1)F)^{\mathbf C^*}\le n+2$.
By the cohomology exact sequence of \eqref{ses1} with $m=n-2$, we obtain an exact sequence
\begin{align}
0\lra H^0((n-3)F)\lra H^0((n-2)F)\lra H^0((n-2)K_S^{-1}).
\end{align}
Since the divisor $S\in |F|$ is $\mathbf C^*$-invariant, this sequence is $\mathbf C^*$-equivariant.
Taking $\mathbf C^*$-fixed part, we obtain an exact sequence
\begin{align}\label{ses4}
0\lra H^0((n-3)F)^{\mathbf C^*}\lra H^0((n-2)F)^{\mathbf C^*}\lra H^0((n-2)K_S^{-1})^{\mathbf C^*}.
\end{align}
Now as is already remarked we have $V_m=H^0(mF)=H^0(mF)^{\mathbf C^*}$ for $m\le n-3$ and $\dim V_m=m+1$ for any $m\ge 0$.
On the other hand, by Prop.\,\ref{prop-multan1} (ii), we have $\dim H^0((n-2)K_S^{-1})^{\mathbf C^*}=1$.
Hence by the exact sequence \eqref{ses4} we obtain $\dim H^0((n-2)F)^{\mathbf C^*}\le (n-1)$.
Taking $m=n-1$ in \eqref{ses1}, the same argument implies an exact sequence
\begin{align}\label{ses5}
0\lra H^0((n-2)F)^{\mathbf C^*}\lra H^0((n-1)F)^{\mathbf C^*}\lra H^0((n-1)K_S^{-1})^{\mathbf C^*}.
\end{align}
Further $\dim H^0((n-2)K_S^{-1})^{\mathbf C^*}=3$ by Prop.\,\ref{prop-multan2} (i).
Hence we obtain $\dim H^0((n-1)F)^{\mathbf C^*}\le (n-1)+3=n+2$.

To finish a proof of (i) and (ii), it remains to see that the $(n+2)$ divisors \eqref{trivial1} and \eqref{nt2} are linearly independent.
The $n$ divisors \eqref{trivial1} are clearly linearly independent and its base locus is exactly the cycle $C$.
On the other hand a ($\mathbf C^*$-invariant) section of the line bundle $(m-1)F$ which defines 
the first one of \eqref{nt2} 
does not belong to $V_{n-1}$, since the divisor does not contain the component $\ol{C}_{n+1}$.
Thus \eqref{trivial1} plus the first divisor of \eqref{nt2} are linearly independent.
It is readily seen that the base locus of this $n$-dimensional subsystem (of $|(n-1)F|^{\mathbf C^*}$) is $C-\ol{C}_{n+1}$.
In particular, it contains $C_{n+1}$.
 On the other hand, the second divisor of \eqref{nt2} does not contain $C_{n+1}$.
 Hence this divisor is not a member of the above $n$-dimensional subsystem.
Thus we have shown that the $(n+2)$ divisors \eqref{trivial1} and \eqref{nt2} are linearly independent. Hence we obtain (i) and (ii).

(iii) is immediate if we note that Bs\,$|V_{n-1}|=C$ and that we can explicitly write the restrictions $S_i^{\pm}|_S$ $(1\le i\le n+1)$ as a sum of the components $C_j$ and $\ol C_j$.
\proofend

\bigskip
We note that as we have $\dim H^0((n-1)F)^{\mathbf C^*}=n+1$ the sequence
\begin{align}\label{ses55}
0\lra H^0((n-2)F)^{\mathbf C^*}\lra H^0((n-1)F)^{\mathbf C^*}\lra H^0((n-1)K_S^{-1})^{\mathbf C^*}\lra 0
\end{align}
is exact.

\subsection{Defining equations of the minitwistor spaces and their singularities}\label{ss-mt}
Thanks to Prop.\,\ref{prop-generator1}, we can determine the image of the meromorphic map $\Phi_{n-1}^{\mathbf C^*}$ (associated to the system $|(n-1)F|^{\mathbf C^*}$) as follows.
(Recall that by the diagram \eqref{cd2}, $\Phi_{n-1}^{\mathbf C^*}(Z)\subset(\pi_{n-1}^{\mathbf C^*})^{-1}(\Lambda_{n-1})$ holds.)

\begin{thm}\label{thm-mt1}
Let $Z$ be as in Prop.\,\ref{prop-fundamental}.
Then there exist a homogeneous coordinate $(z_1,z_2,\cdots,z_{n})$ on $\mathbf P^{\vee}V_{n-1}=\mathbf{CP}^{n-1}$ and two sections $z_{n+1},z_{n+2}\in H^0((n-1)F)^{\mathbf C^*}$ which satisfy the following.
(i) $\{z_i\set 1\le i\le n+2\}$ is a basis of $H^0((n-1)F)^{\mathbf C^*}$.
(ii) With respect to the non-homogeneous coordinate $(z_2/z_1,z_3/z_1,\cdots,z_n/z_1)$ on \ $\mathbf P^{\vee}V_{n-1}$, the rational normal curve $\Lambda_{n-1}$ in $\mathbf P^{\vee}V_{n-1}$ is given by
\begin{align}\label{rnc1}
\left\{\left(\lambda,\lambda^2,\cdots,\lambda^{n-1}\right)\set\lambda\in\mathbf C\right\}.
\end{align} 
(iii) $z_{n+1}$ and $z_{n+2}$ define the 2 divisors \eqref{nt2} respectively and also satisfy $\ol{z}_{n+2}=z_{n+1}$.
(iv) 
The image $\mathscr T:=\Phi_{n-1}^{\mathbf C^*}(Z)$ is a surface which satisfies not only  relations in the defining ideal of $\Lambda_{n-1}$, but also the following quadratic equation
\begin{align}\label{nt3}
z_{n+1}z_{n+2}=z_2\left(
z_n-\sigma_1z_{n-1}+\sigma_2z_{n-2}+\cdots+(-1)^{n-1}\sigma_{n-1}z_1
\right),
\end{align}
where $\sigma_i$ is the elementary symmetric polynomial  of degree $i$ of some $(n-1)$ real numbers $\lambda_3,\cdots,\lambda_{n+1}$.
(v) The degree of the surface $\mathscr T$ in $\mathbf{CP}^{n+1}$ is $2(n-1)$.
\end{thm}

\noindent
Proof.
Let $e_i$ ($1\le i\le n+1$) be a section of the line bundle $[S_i^+]$ which defines $S_i^+$.
We put
\begin{align}
y_i=e_i\cdot \ol{e}_i\in H^0(F)\hspace{5mm} (1\le i\le n+1).
\end{align}
(Of course, this defines the divisor $S_i^++S_i^-$.)
We choose $\{y_1,y_2\}$ as a basis of $H^0(F)$.
Since these are also a basis of the real part $H^0(F)^{\sigma}$, there exist real numbers $\lambda_3,\cdots,\lambda_{n+1}$ such that 
\begin{align}
y_i=y_2-\lambda_iy_1
\end{align}
hold.
We put $\lambda_2=0$.
Then since $S_i^++S_i^-\in |F|^{\sigma}\simeq S^1$ are arranged in a linear order, either $\lambda_2=0<\lambda_3<\lambda_4<\cdots<\lambda_{n+1}$ or $\lambda_2=0>\lambda_3>\lambda_4>\cdots>\lambda_{n+1}$ holds.
We put 
\begin{align}\label{z1}
z_1=y_1^{n-1},\,z_2=y_1^{n-2}y_2,\,\cdots\,,z_n=y_2^{n-1}.
\end{align}
Then it is obvious from the definition of the space $V_{n-1}$ that $\{z_1,z_2,\cdots,z_{n}\}$ is a basis of $V_{n-1}$.
Also it is immediate to see that the curve $\Lambda_{n-1}$ is represented as in \eqref{rnc1}.
We put 
\begin{align}
z_{n+1}=e_1^{n-2}e_2e_3\cdots e_{n+1}\,\,\,\text{and}\,\,\,\, z_{n+2}=\ol{e}_1^{n-2}\ol{e}_2\ol{e}_3\cdots \ol{e}_{n+1}.
\end{align}
Then by Prop.\,\ref{prop-generator1}, $\{z_1,\cdots,z_{n+2}\}$ is a basis of $H^0((n-1)F)^{\mathbf C^*}$. 
Also it is obvious that $z_{n+1}$ and $z_{n+2}$ define the divisors \eqref{nt2} respectively.
Then we have
\begin{align}
z_{n+1}z_{n+2}&= (e_1\ol{e}_1)^{n-2}\cdot (e_2\ol{e}_2)\cdot (e_3\ol{e}_3)\cdots(e_{n+1}\ol{e}_{n+1})\\
&=y_1^{n-2}y_2y_3\cdots y_{n+1}\\
&=y_1^{n-2}y_2\cdot(y_2-\lambda_3y_1)(y_2-\lambda_4y_1)\cdots(y_2-\lambda_{n+1}y_1).
\label{dcp1}
\end{align} 
Expanding the right-hand side and using \eqref{z1}, we obtain the equation \eqref{nt3}.
Thus we obtain (iv).
For (v) note that the surface $\mathscr T$ is an intersection of a rational scroll of lines whose degree is $(n-1)$ (= the degree of $\Lambda_{n-1}$ in $\mathbf{CP}^{n-1}$) and the hyperquadric \eqref{nt3}.
The claim is immediate from this.
%
\proofend

\bigskip
Since the restriction map in the sequence \eqref{ses55} is surjective, Prop.\,\ref{prop-multan2} (v) implies the following

\begin{prop}\label{prop-quot1}
Restriction of the meromorphic map $\Phi_{n-1}^{\mathbf C^*}$ to a general member  $S\in |F|$ is a morphism and its general  fibers are irreducible smooth rational curves in $S$.
(Namely the restriction can be regarded as a (holomorphic) quotient map of $\mathbf C^*$-action on $S$.)
\end{prop}

By the diagram \eqref{cd2}, we have the following commutative diagram of meromorphic maps
\begin{equation}\label{cd3}
 \CD
Z@>{\Phi_{n-1}^{\mathbf C^*}}>>\mathscr T \\
 @V\Psi_{n-1}VV @VV{\pi_{n-1}^{\mathbf C^*}}V\\
\Lambda_{n-1}@>>>\mathbf{CP}^{n-1}\\
 \endCD
 \end{equation}
 By Prop.\,\ref{prop-quot1}, the map $\Phi_{n-1}^{\mathbf C^*}:Z\to\mathscr T$ can be regarded as a quotient map by $\mathbf C^*$-action.
 Namely the complex surface $\mathscr T$ can be considered as a minitwistor space of the present twistor space $Z$.
 
 Next we investigate geometric structure of the surface $\mathscr T$.
 The equation \eqref{nt3} means that  by the projection $\pi_{n-1}^{\mathbf C^*}$, $\mathscr T$ is birational to a conic bundle over $\Lambda_{n-1}$ whose discriminant locus is  the intersection of $\Lambda_{n-1}$ and the union of 2 hyperplanes determined as the zero locus of the right-hand side of \eqref{nt3}.
Since $z_2=y_1^{n-2}y_2$ as in \eqref{z1}, the rational normal curve $\Lambda_{n-1}$ intersects the hyperplane $\{z_2=0\}$ transversally at $(1,0,0,\cdots,0)$, and  touches the hyperplane at $(0,0,\cdots,0,1)$ with multiplicity $(n-2)$, 
where we are using the homogeneous coordinate $(z_1,\cdots,z_n)$ on $\mathbf{CP}^{n-1}$ given in Theorem\,\ref{thm-mt1}.
 Also, by \eqref{dcp1}, the intersection of $\Lambda_{n-1}$ with another hyperplane consists of $(n-2)$ points $(\lambda_i,\lambda_i^2,\cdots,\lambda_i^{n-1})$ ($3\le i\le n+1$) with respect to the  coordinate on $\Lambda_{n-1}$ in \eqref{rnc1}.
The fibers of these $(n+1)$ discriminant points consist of two irreducible components.
 The reducible fiber over the point $(\lambda_i,\lambda_i^2,\cdots,\lambda_{n+1}^2)\in\Lambda_{n-1}$ is precisely $\Phi(S_i^+)+\Phi(S_i^-)$ for $2\le i\le n+1$ and that over the touching point $(0,\cdots,0,1)\in\Lambda_{n-1}$ is $\Phi(S_1^+)+\Phi(S_1^-)$.
(See also Figure \ref{fig-mt}.)

Note that the complex surface $\mathscr T$ is uniquely determined by the set $\{\lambda_3,\cdots,\lambda_{n+1}\}$ of different $n-1$ numbers.
Further, two surfaces $\mathscr T$ and $\mathscr T'$ determined by $\{\lambda_3,\cdots,\lambda_{n+1}\}$  and $\{\lambda'_3,\cdots,\lambda'_{n+1}\}$ respectively are isomorphic iff there exists a constant $c\in\mathbf R^{\times}$ such that $\{c\lambda_3,\cdots,c\lambda_{n+1}\}=\{\lambda'_3,\cdots,\lambda'_{n+1}\}$.
Thus the present minitwistor spaces form $(n-2)$-dimensional family.

It is obvious that with respect to the homogeneous coordinates  $(z_1,\cdots,z_n)$ on $\mathbf P^{\vee}V_{n-1}$ and $(z_1,\cdots,z_{n+2})$ on $\mathbf P^{\vee} H^0((n-1)F)^{\mathbf C^*}$, the projection $\pi_{n-1}^{\mathbf C^*}$ is explicitly given by
$\pi_{n-1}^{\mathbf C^*}(z_1,\cdots,z_{n+2})=(z_1,\cdots,z_n).$
Hence the center of $\pi_{n-1}^{\mathbf C^*}$ is the line defined by 
\begin{align}
l_{\infty}:=\{z_1=z_2=\cdots=z_n=0\}.
\end{align}
Fibers of $\pi_{n-1}^{\mathbf C^*}$ are planes containing $l_{\infty}$.
It can be readily computed by Theorem\,\ref{thm-mt1} that we have
\begin{align}\label{inters1}
l_{\infty}\cap \mathscr T=\{(0,0,\cdots,0,1,0), (0,0,\cdots,0,0,1)\}
\end{align}
in the homogeneous coordinate $(z_1,\cdots,z_{n+2})$.
Fibers of $\pi_{n-1}^{\mathbf C^*}|_{\mathscr T}$ form a pencil on $\mathscr T$ whose inverse image by $\Phi_{n-1}^{\mathbf C^*}$ are the pencil $|F|$ on $Z$.
The base locus of this pencil on $\mathscr T$ is precisely the 2 points \eqref{inters1}.

Singularities of the image surface $\mathscr T$ is described as follows.

\begin{prop}\label{prop-mt2}
Let $\mathscr T$ be as in Theorem\,\ref{thm-mt1}.
Then the singularities of $\mathscr T$ consist of the following 3 points.
(a) the 2 intersection points $l_{\infty}\cap \mathscr T$ in \eqref{inters1},
(b) the point $(0,\cdots,0,1,0,0)$ in the homogeneous coordinate $(z_1,\cdots,z_{n+2})$.
Further, the 2 points of (a) are isomorphic to the cyclic quotient singularity of $\mathbf Z_{n-1}$-action $(z,w)\mapsto (\zeta z,\zeta w)$ on $\mathbf C^2$ with $\zeta=e^{\frac{2\pi\sqrt{-1}}{n-1}}$.
The point (b) is isomorphic to a rational double point of A$_{n-3}$-type.
\end{prop}

We omit a proof since the description of $\mathscr T$ obtained in Theorem\,\ref{thm-mt1} is completely explicit and by using it, it is not difficult to derive the conclusion.
Instead we only note that the fiber of $\pi_{n-1}^{\mathbf C^*}:\mathscr T\to\mathbf{CP}^{n-1}$ in the diagram \eqref{cd3}  over the (touching) point $(z_1,\cdots,z_n)=(0,\cdots,0,1)$ consists of 2 irreducible component, and that the point in (b) is precisely the intersection point of these components. 

Let $\hat{\mathscr T}\to\mathscr T$ be the blowing-up of the embedded surface $\mathscr T\subset\mathbf{CP}^{n+1}$ along the line $l_{\infty}$.
Then the pair of  cyclic quotient singularities ((a) of Prop.\,\ref{prop-mt2}) are resolved.
At the same time, the indeterminacy locus of the projection $\mathscr T\to\Lambda_{n-1}$ is resolved so that $\hat{\mathscr T}$ has a structure of a conic bundle over $\Lambda_{n-1}$.
Let $\Gamma$ and $\ol{\Gamma}$ be the exceptional curves over the 2 points \eqref{inters1}.
These are of course  smooth rational curves in $\hat{\mathscr T}$ whose self-intersection numbers  are $(1-n)$.
$\hat{\mathscr T}$ still has $A_{n-3}$-singularity ((b) of Prop.\,\ref{prop-mt2}).
Let $\tilde{\mathscr T}\to\hat{\mathscr T}$ be its minimal resolution.
The exceptional curve of this resolution is a string of $(-2)$-curves consisting of $(n-3)$ irreducible components.
The composition
\begin{align}\label{resol1}
\nu:\tilde{\mathscr T}\lra\hat{\mathscr T}\lra\mathscr T
\end{align}
 is the minimal resolution of the surface $\mathscr T$.
Since the natural projection $\tilde{\mathscr T}\to\Lambda_{n-1}$ has precisely $(n+1)$ reducible fibers, and since the number of their irreducible components is $2n+(n-1)=3n-1=(2n-2)+(n+1)$, we obtain
\begin{align}\label{chern}
c_1^2(\tilde{\mathscr T})=8-(2n-2)=10-2n.
\end{align}
From this we obtain $b_2(\tilde{\mathscr T})=2n$.

\section{Description of the twistor spaces as conic bundles}
\subsection{Partial elimination of the indeterminacy locus of the quotient meromorphic map}\label{ss-pr}
We recall that the surface $\mathscr T$ studied in \S \ref{ss-mt} was the image of the meromorphic map $\Phi_{n-1}^{\mathbf C^*}$ associated to the linear system $|(n-1)F|^{\mathbf C^*}$, and it can also be regarded as a quotient space of $\mathbf C^*$-action on the twistor space $Z$.
 Since we know generator of this system  in an explicit form as in 
Prop.\,\ref{prop-generator1}, we can in principle eliminate the indeterminacy locus of the map $\Phi_{n-1}^{\mathbf C^*}$, by a succession of blowing-ups. 
As in Prop.\,\ref{prop-generator1} we have Bs\,$|(n-1)F|^{\mathbf C^*}=C-C_{n+1}-\ol{C}_{n+1}$.
Because the generators in the proposition intersect in a quite complicated way,
 it looks not easy to give a complete elimination.
So from now on we give an elimination only in a neighborhood of the curve $C_2\cup\ol{C}_2$.
This is enough for the purpose of obtaining  projective models of the present twistor spaces as conic bundles.
We recall that $C_2$ and $\ol{C}_2$ are only  components among the cycle $C$ which are pointwisely fixed by the $\mathbf C^*$-action.
In the terminology of \cite{S82}, $C_2$ and $\ol{C}_2$ are `source' and `sink' of the $\mathbf C^*$-action on $Z$.
(Namely, by the $\mathbf C^*$-action, general points on $Z$ goes to points on $C_2$ and $\ol{C}_2$ as $s\to 0$ and $s\to\infty$.)
In the following we repeat the next usual way of  eliminations of the base locus:
\begin{itemize}
\item blow-up along base curves,
\item compute the total transforms of generators,
\item compute their fixed components
\item remove the fixed component from the total transforms
\item  compute the base curves
\end{itemize}

As the first step, let $Z_1\to Z$ be blowing-up along $C_2\cup \ol{C}_2$, and $E_2$ and $\ol{E}_2$ the exceptional divisors over $C_2$ and $\ol{C}_2$ respectively.
Since the 2 divisors $S_2^+$ and $S_3^-$ (in $Z$) intersect transversally along $C_2$, and the self-intersection numbers of $C_2$ in  $S_2^+$ and $S_3^-$ are $(-1)$ for both,
the normal bundle of $C_2$ in $Z$ is isomorphic to $\mathscr O(-1)^{\oplus 2}$.
Hence $E_2$  is biholomorphic to the trivial $\mathbf{CP}^1$-bundle over $C_2$.
Similarly, $\ol{E}_2\simeq \ol{C}_2\times\mathbf{CP}^1$.
Then since every generators in Prop.\,\ref{prop-generator1} contain the curve $C_2\cup\ol{C}_2$ with the same multiplicity $(n-1)$, 
their total transforms in $Z_1$ contain $E_2+\ol{E}_2$ with multiplicity $(n-1)$.
Hence subtracting  $(n-1)E_2+(n-1)\ol{E}_2$ from the total transforms of generators,
we obtain a linear system on $Z_1$ 
whose generators can still be written as \eqref{trivial1} and \eqref{nt2}, where
 we are using the same notations to denote the divisors in $Z$ and their strict transforms in $Z_1$ as promised in the convention.
 
 This $(n+1)$-dimensional linear system on $Z_1$ still has a base locus which intersects $E_2$ or $\ol{E}_2$.
 Namely, the curves $C_1$ and $C_3$ (in $Z_1$), intersecting $E_2$ transversally at a unique point for each, are (still)  base curves.
 Similarly, $\ol{C}_1$ and $\ol{C}_3$ (in $Z_1$) are base curves intersecting $\ol{E}_2$ transversally at a unique point  for each.
 More precisely, all divisors in \eqref{trivial1} (viewed as divisors on $Z_1)$ obviously contain $C_1,\ol{C}_1,C_3$ and $\ol{C}_3$ with multiplicity $(n-1)$.
 On the other hand, the former divisor of \eqref{nt2} contains these 4 curves with multiplicity $(n-2)$, $n$, $n$ and $(n-2)$ respectively.
 The latter divisor of \eqref{nt2} contains these curves with multiplicities $n,n-2,n-2,n$ respectively.
 It can also be  verified (by using Prop.\,\ref{prop-generator1}) that there is no base locus other than these which intersects $E_2$ or $\ol{E}_2$.
  Let $Z_2\to Z_1$ be the blowing-up along $C_1\cup \ol{C}_1\cup C_3\cup\ol{C}_3$, and $E_1, \ol{E}_1, E_3$ and $\ol{E}_3$ their exceptional divisors respectively.
  Then  all the total transforms of the divisors \eqref{trivial1} and \eqref{nt2} into $Z_1$ contain $E_1, \ol{E}_1, E_3 $ and $\ol{E}_3$  with the above multiplicities respectively.
  So removing the maximal common divisor $(n-2)(E_1+E_3+\ol{E}_1+\ol{E}_3)$
  (namely, the fixed component) from the total transforms, we obtain a linear system on $Z_2$ (with the same dimension $(n+1)$) whose generators are 
  \begin{align}\label{trivial2}
(n-1-k)(S_1^++S_1^-)+k(S_2^++S_2^-)+(E_1+E_3+\ol{E}_1+\ol{E}_3),\,\,\,0\le k\le n-1,
\end{align}
\begin{align}\label{nt4}
(n-2)S_1^+\,+\sum_{i=2}^{n+1} S_i^++2E_3+2\ol{E}_1,\,\,\,
(n-2)S_1^-\,+\sum_{i=2}^{n+1} S_i^-+2E_1+2\ol{E}_3.
\end{align}
Of course, \eqref{trivial2} and \eqref{nt4}  correspond \eqref{trivial1} and \eqref{nt2} respectively.
From these, we can deduce that the base locus  which intersect $E_2$ are the following $(n-1)$ rational curves.
\begin{align}\label{base1}
S_1^+\cap E_1\,\,\,\,{\text{and}}\,\,\,S_i^-\cap E_3\,\,(4\le i\le n+1).
\end{align}
Similarly, the base locus of the linear system (on $Z_2$) intersecting $\ol{E}_2$ are the following $(n-1)$ rational curves
\begin{align}\label{base2}
S_1^-\cap \ol{E}_1\,\,\,\,{\text{and}}\,\,\,S_i^+\cap \ol{E}_3\,\,(4\le i\le n+1).
\end{align}
Let $Z_3\to Z_2$ be the blowing-up along $2(n-1)$ rational curves \eqref{base1} and \eqref{base2}, $F_1$ and $D_i$ ($4\le i\le n+1$) the exceptional divisors over the curves \eqref{base1} respectively, and $\ol{F}_1$ and $\ol{D}_i$  ($i\le 4\le n+1$) the exceptional divisors over the curves \eqref{base2} respectively.
Counting the multiplicity of the generators \eqref{trivial2} and \eqref{nt4} along the base curves \eqref{base1} and \eqref{base2} respectively, we deduce that the fixed component of the $(n+1)$-dimensional linear system on $Z_3$ is precisely
 $(F_1+\ol{F}_1)+\sum_{i=4}^{n+1}(D_i+\ol{D}_i)$.
 Removing this  from the total transforms of the generators \eqref{trivial2} and \eqref{nt4}, 
as  generators of the movable part of  the linear system on $Z_3$, we can choose the following divisors:
   \begin{align}\label{trivial3}
(n-1-k)(S_1^++S_1^-)+k(S_2^++S_2^-)+(E_1+E_3+\ol{E}_1+\ol{E}_3)+(n-1-k)(F_1+\ol{F}_1),
\end{align}
with $0\le k\le n-1$, and 
\begin{align}\label{nt5}
(n-2)S_1^+\,+\sum_{i=2}^{n+1} S_i^++2E_3+2\ol{E}_1+(n-3)F_1+\sum_{i=4}^{n+1}D_i+\ol{F}_1,
\end{align}
\begin{align}\label{nt6}
(n-2)S_1^-\,+\sum_{i=2}^{n+1} S_i^-+2\ol{E}_3+2E_1+F_1+(n-3)\ol{F}_1+\sum_{i=4}^{n+1}\ol{D}_i.
\end{align}
This linear system on $Z_3$ still has base curves intersecting $E_2\cup\ol E_2$. 
Namely, the  two rational curves $E_1\cap F_1$ and $\ol{E}_1\cap \ol{F}_1$ are such curves.
So let $Z_4\to Z_3$ be the blowing-up along these 2 curves, and $F_2$ and $\ol{F}_2$ the exceptional divisors over there respectively.
Pulling back the divisors \eqref{trivial3}--\eqref{nt6} and removing its fixed component $F_2+\ol{F}_2$, we obtain, as  generators of an $(n+1)$-dimensional linear system on $Z_4$, 
the following divisors:
  \begin{align}\label{trivial4}
(n-1-k)(S_1^++S_1^-)+k(S_2^++S_2^-)+(E_1+E_3+\ol{E}_1+\ol{E}_3)+(n-1-k)\sum_{i=1}^2(F_i+\ol{F}_i)
\end{align}
with $0\le k\le n-1$, and 
\begin{align}\label{nt7}
(n-2)S_1^+\,+\sum_{i=2}^{n+1} S_i^++2E_3+2\ol{E}_1+(n-3)F_1+(n-4)F_2+\sum_{i=4}^{n+1}D_i+\ol{F}_1,
\end{align}
\begin{align}\label{nt8}
(n-2)S_1^-\,+\sum_{i=2}^{n+1} S_i^-+2\ol{E}_3+2E_1+F_1+(n-3)\ol{F}_1+(n-4)\ol{F}_2+\sum_{i=4}^{n+1}\ol{D}_i.
\end{align}
If $n=4$, since the coefficients of $F_2$ in   \eqref{nt7} and $\ol{F}_2$ in \eqref{nt8} are both zero,  this system on $Z_4$ is free, at least in a neighborhood of $E_2\cup\ol{E}_2$.
So we stop the operations here.
If $n\ge 5$, the base locus of this system which intersects  $E_2\cup\ol{E}_2$ consists of 2 rational curves $E_1\cap F_2$ and $\ol{E}_1\cap\ol{F}_2$.
So let $Z_5\to Z_4$ be the blowing-up along these 2 curves with the exceptional divisor $F_3$ and $\ol{F}_3$ respectively.
When $n=5$, the linear system on $Z_5$ obtained by pulling back the divisors  \eqref{trivial4}--\eqref{nt8} and removing the fixed component  $F_3+\ol{F}_3$ becomes free, and we stop the operations here.
For $n\ge 6$,  repeating this operation $(n-4)$ times,  we obtain a string of 3-folds $Z_{n}\to Z_{n-1}\to\cdots \to Z_5$ and exceptional  divisors $F_4,\ol{F}_4,F_5,\ol{F}_5,\cdots F_{n-2},\ol{F}_{n-2}$, together with an $(n+1)$-dimensional linear system on $Z_n$ whose generators are given by
  \begin{align}\label{trivial5}
(n-1-k)(S_1^++S_1^-)+k(S_2^++S_2^-)+(E_1+E_3+\ol{E}_1+\ol{E}_3)+(n-1-k)\sum_{i=1}^{n-2}(F_i+\ol{F}_i)
\end{align}
with $0\le k\le n-1$, and 
\begin{align}\label{nt9}
(n-2)S_1^+\,+\sum_{i=2}^{n+1} S_i^++2E_3+2\ol{E}_1+\sum_{i=1}^{n-2}(n-2-i)F_i+\sum_{i=4}^{n+1}D_i+\ol{F}_1,
\end{align}
\begin{align}\label{nt10}
(n-2)S_1^-\,+\sum_{i=2}^{n+1} S_i^-+2\ol{E}_3+2E_1+F_1+\sum_{i=1}^{n-2}(n-2-i)\ol{F}_i+\sum_{i=4}^{n+1}\ol{D}_i.
\end{align}
Then basically because the coefficients of $F_{n-2}$ in \eqref{nt9} and  $\ol{F}_{n-2}$ in \eqref{nt10} 
are both zero, this linear system on $Z_n$ is free, at least  in a neighborhood of the divisor $E_2\cup\ol{E}_2$.
Thus we have obtained an explicit sequence of blowing-ups $Z_n\to Z_{n-1}\to\cdots\to Z_1\to Z$ which eliminates the base locus of the system $|(n-1)F|^{\mathbf C^*}$, at least in a neighborhood of $C_2\cup\ol{C}_2$.

Then we have the following

\begin{prop}\label{prop-elim1}
Consider the following diagram 
\begin{equation}\label{cd4}
 \CD
Z_n@>>> Z \\
 @V\tilde{\Phi}_{n-1}^{\mathbf C^*}VV  @VV{\Phi_{n-1}^{\mathbf C^*}}V\\
\tilde{\mathscr T}@>{\nu}>>\mathscr T\\
 \endCD
 \end{equation}
 where $\nu$ is the minimal resolution of $\mathscr T$ as in \eqref{resol1}, and 
$\tilde{\Phi}_{n-1}^{\mathbf C^*}$ is a  meromorphic map  uniquely determined by the commutativity of the diagram.
Then  $\tilde{\Phi}_{n-1}^{\mathbf C^*}$ is holomorphic in  neighborhoods of $E_2$ and $\ol{E}_2$, and its restrictions onto these divisors are both biholomorphic.
\end{prop}

\noindent
Proof.
We consider the pencil $|F|$ on $Z$ whose base locus is the cycle $C$ as in Prop.\,\ref{prop-fundamental} (i).
Since the first operation $Z_1\to Z$ blow-ups $C_2\cup\ol{C}_2\subset C$,
the strict transform of the pencil $|F|$ into $Z_1$ makes sense.
Similarly, since $Z_2\to Z_1$ blow-ups $C_1\cup\ol{C}_1\cup C_3\cup\ol{C}_3$ which are still contained in the base locus of the pencil on $Z_1$,
the strict transform of the pencil into $Z_2$ also makes sense.
Evidently this pencil on $Z_2$ has no base point, at least in a neighborhood of $E_2\cup\ol{E}_2$.
Pulling this back by the sequence of blowing-ups $Z_n\to Z_{n-1}\to\cdots\to Z_2$, we obtain a pencil on $Z_n$ which has no base point in a neighborhood of $E_2\cup\ol{E}_2$.
By restricting this pencil to $E_2$ and $\ol{E}_2$ ($\subset Z_n$), we obtain pencils on 
 $E_2$ and $\ol{E}_2$ without base points.
 
 Since the sequence  $Z_n\to Z_{n-1}\to\cdots\to Z$ eliminates the base locus of the system $|(n-1)F|^{\mathbf C^*}$ in a neighborhood of $C_2\cup \ol{C}_2$, the composition map $Z_n\to Z\to\mathscr T$ is holomorphic in a neighborhood of $E_2\cup\ol{E}_2$.
 By the commutative diagram \eqref{cd3}, members of the above pencils on $E_2$ and $\ol{E}_2$ are mapped to fibers of the (rational) conic bundle $\pi_{n-1}^{\mathbf C^*}:\mathscr T\to\Lambda_{n-1}$.
 General (irreducible) members of the pencils are mapped biholomorphically to smooth fiber of $\mathscr T\to\Lambda_{n-1}$.
Further, from our explicit way of the sequence of blow-ups, it is obvious that these pencils (on $E_2$ and $\ol{E}_2$) have precisely $(n+1)$ reducible members, each of which are explicitly given by the restriction of the divisors
 \begin{align}
 S_1^++\sum_{i=1}^{n-2}F_i,\,\,S_2^++S_2^-,\,\,S_3^++S_3^-,\,\,S_i^-+\ol{D}_i\,\,(4\le i\le n+1)
 \end{align}
 to $E_2$, and the restrictions of their conjugates to $\ol{E}_2$.
 On the other hand, as is already explained in \S\ref{ss-mt},
 reducible fibers of the (rational) conic bundle  $\pi_{n-1}^{\mathbf C^*}:\mathscr T\to\Lambda_{n-1}$ are over the points $(\lambda_i,\lambda_i^2,\cdots,\lambda_i^{n-1})\in\Lambda_{n-1}$, $2\le i\le n+1$, and the touching point $(0,\cdots,0,1)\in\Lambda_{n-1}$.
 Further the fiber over $(\lambda_i,\lambda_i^2,\cdots,\lambda_i^{n-1})$ is $\Phi_{n-1}^{\mathbf C^*}(S_i^+)+\Phi_{n-1}^{\mathbf C^*}(S_i^-)$ for $2\le i\le n+1$ and that over $(0,\cdots,0,1)$ is 
 $\Phi_{n-1}^{\mathbf C^*}(S_1^+)+\Phi_{n-1}^{\mathbf C^*}(S_1^-)$.
 Therefore, the map $E_2\to\mathscr T$ (which is the restriction of the composition $Z_n\to Z\to\mathscr T$) is the birational morphism which contracts the following rational curves:
 \begin{align}
 E_2\cap E_1,\,\,E_2\cap E_3\,\,\,\text{and}\,\,\,F_i\cap E_2\,\,\,\text{for}\,\,\,1\le i\le n-3.
  \end{align}
 By the configuration of these curves (which can be read off from our explicit way of the sequence $Z_n\to\cdots Z)$, it follows that the morphism $E_2\to\tilde{\mathscr T}$ is precisely the minimal resolution of $\mathscr T$.
 Hence we obtain $E_2\simeq\tilde{\mathscr T}$, by the uniqueness of the minimal resolution.
 By reality, we also obtain $\ol E_2\simeq\tilde{\mathscr T}$.
 
 Finally we show that $\tilde{\Phi}_{n-1}^{\mathbf C^*}:Z_n\to\tilde{\mathscr T}$ in the proposition is holomorphic  in neighborhoods of $E_2$ and $\ol{E}_2$.
 By our choice of the sequence $Z_n\to Z_{n-1}\to\cdots Z_1\to Z$, only $E_1,E_3,\ol{E}_1$ and $\ol{E}_2$ are mapped to the conjugate pair of singular points of $\mathscr T$ ((a) of Prop.\,\ref{prop-mt2}), and only $F_i$ and $\ol{F}_i$ $(1\le i\le n-3)$ are mapped to the $A_{n-3}$-singularity of $\mathscr T$ ((b) of Prop.\,\ref{prop-mt2}), by the composition $Z_n\to Z\to\mathscr T$.
 Further, all intersections of these exceptional divisors with $E_2$ and $\ol{E}_2$ are precisely the exceptional curves of the minimal resolutions $E_2\to\mathscr T$ or $\ol{E}_2\to\mathscr T$.
 These imply that the map $Z_n\to Z\to \mathscr T$ is lifted to a map $Z_n\to \tilde{\mathscr T}$ in a way that it is still holomorphic on $E_2$ and $\ol{E}_2$.
%
%
\proofend

\subsection{Description of the discriminant locus of the quotient map}\label{ss-dl}
Next based on Prop.\,\ref{prop-elim1} 
we give projective models of our twistor spaces as conic bundles over the minimal resolution of the minitwistor spaces.
Namely we explicitly construct a $\mathbf{CP}^2$-bundle over the minimal resolution $\tilde{\mathscr T}$, and show that our twistor space $Z$ is bimeromorphically embedded into this bundle as a conic bundle.
We also give a defining equation of the conic bundle.

Although the map $\tilde{\Phi}_{n-1}^{\mathbf C^*}:Z_n\to\tilde{\mathscr T}$ in Prop.\,\ref{prop-elim1} is holomorphic on a neighborhood of $E_2\cup\ol{E}_2$, 
it  still has indeterminacy locus because there are still a base locus.
In fact, the curve $C_4\cup\ol{C}_4$ in $Z_n$ is still a base curve and it is not resolved even after blowing-up along $C_4\cup\ol{C}_4$.
(Namely another base locus appears on the exceptional divisors over $C_4$ and $\ol{C}_4$.)
Also, there remain many base curves on the exceptional divisors $E_3$ and $\ol{E}_3$ (in $Z_n$) which have a complicated structure.
To remove these base curves completely, we need a lot of blow-ups and it looks difficult to give them in an explicit form.
We do not persist in them and  take any  sequence of blowing-ups $\tilde Z\to Z_n$ along $\mathbf C^*$-invariant non-singular centers which eliminates the  base locus of the system.
We can suppose that all centers of the blow-ups are disjoint from $E_2\cup\ol{E}_2$ since $\tilde{\Phi}_{n-1}^{\mathbf C^*}$ is already holomorphic on these divisors.
We denote the resulting morphism by $\tilde{\Phi}^{\mathbf C^*}:\tilde{Z}\to\tilde{\mathscr T}$.
General fibers of $\tilde{\Phi}^{\mathbf C^*}$ are $\mathbf C^*$-invariant irreducible rational curves by Prop.\,\ref{prop-quot1}, and  $E_2$ and $\ol{E}_2$ are sections of $\tilde{\Phi}^{\mathbf C^*}$ which are fixed by the $\mathbf C^*$-action on $\tilde Z$.
We consider the direct image sequence of the exact sequence 
\begin{align}
0\,\lra\, \mathscr O_{\tilde Z}\,\lra\,\mathscr O_{\tilde Z}(E_2+\ol{E}_2)\,\lra\,N_{E_2/\tilde{Z}}\oplus N_{\ol{E}_2/\tilde{Z}}\,\lra\,0.
\end{align}
Since every fiber of $\tilde{\Phi}^{\mathbf C^*}$ is at most a string of rational curves,
we have $R^1\tilde\Phi^{\mathbf C^*}\mathscr O_{\tilde Z}=0$.
Hence a part of the direct image sequence becomes
\begin{align}\label{ses50}
0\,\lra\, \mathscr O_{\tilde{\mathscr  T}}\,\lra\,(\tilde\Phi^{\mathbf C^*})_*\mathscr O_{\tilde Z}(E_2+\ol{E}_2)\,\lra\,N_{E_2/\tilde{Z}}\oplus N_{\ol{E}_2/\tilde{Z}}\,\lra\,0,
\end{align}
where $N_{E_2/\tilde{Z}}$ and $ N_{\ol{E}_2/\tilde{Z}}$ are considered as line bundles over $\tilde{\mathscr T}$.
On the other hand there are obvious isomorphisms $N_{E_2/\tilde{Z}}\simeq N_{E_2/Z_n}$ and $N_{\ol{E}_2/\tilde{Z}}\simeq N_{\ol{E}_2/Z_n}$.
Further since our sequence of blow-ups $Z_n\to Z_{n-1}\to \cdots\to Z_1\to Z$ are explicit, we can concretely compute the normal bundles $N_{E_2/Z_n}$ and $ N_{\ol{E}_2/Z_n}$.
In particular, basically by the reason that the degree of $N_{E_2/Z_n}$ along general fibers of the natural projection $E_2\to C_2$ is $-1$, we obtain 
\begin{align}
H^1(N_{E_2/Z_n})=H^1(N_{\ol{E}_2/Z_n})=0.
\end{align}
These imply that the exact sequence \eqref{ses50} splits and we obtain an isomorphism
\begin{align}
(\tilde\Phi^{\mathbf C^*})_*\mathscr O(E_2+\ol{E}_2)
\simeq
N_{E_2/Z_n}\oplus N_{\ol{E}_2/Z_n}\oplus\mathscr O_{\tilde{\mathscr T}}.
\end{align}
Then let
\begin{equation}\label{p2bdle}
\mu: \tilde{Z}\lra\mathbf{P}( N_{E_2/Z_n}^{\vee}\oplus N_{\ol{E}_2/Z_n}^{\vee}\oplus\mathscr O)
\end{equation}
be the relative meromorphic map over $\tilde{\mathscr T}$ associated to the pair $\{\mathscr O(E_2+\ol{E}_2),\tilde{\Phi}^{\mathbf C^*}\}$.
Obviously $\mu$ is bimeromorphic over its image and
the image $\mu(\tilde Z)$ is a conic bundle (over $\tilde{\mathscr T}$).
The discriminant locus of the projection $\mu(\tilde{Z})\to \tilde{\mathscr T}$ is a member of the system $|N_{E_2/Z_n}^{\vee}\otimes N_{\ol{E}_2/Z_n}^{\vee}|$.
Note that this system can be explicitly determined.
For the purpose of determining the discriminant curves, we show the following.

\begin{lemma}\label{lemma-hps1}
Let $L_0$ be the fixed twistor line in $Z$ (Def.\,\ref{def-ftl}). 
Then we have the following.
(i) $L_0$ is disjoint from the base locus of the system $|(n-1)F|^{\mathbf C^*}$, and the image
$\mathscr C_0:=\Phi_{n-1}^{\mathbf C^*}(L_0)$ is a curve on $\mathscr T$. 
(ii) $\mathscr C_0$ is a hyperplane section of $\mathscr T$ with respect to the natural embedding $\mathscr T\subset\mathbf{CP}^{n+1}$.
(iii) The virtual genus of the curve $\mathscr C_0$ is $n-2$.
\end{lemma}
\noindent
Since $\mathscr C_0$ is a rational curve, (iii) means that $\mathscr C_0$ is a singular curve.

\noindent
Proof of Lemma \ref{lemma-hps1}.
By Prop.\,\ref{prop-generator1}, we have Bs\,$|(n-1)F|^{\mathbf C^*}=C-C_{n+1}-\ol{C}_{n+1}$.
$L_0$ is not contained in this base curve since $L_0$ is itself real.
If $L_0$ intersects this base curve, the intersection must be a $\mathbf C^*$-fixed point.
However, twistor lines going through such a point  on $C$ are not $\mathbf C^*$-fixed, as in Prop.\,\ref{prop-u1}.
Hence $L_0\cap$\,Bs\,$|(n-1)F|^{\mathbf C^*}=\emptyset$.
Then the image $\mathscr C_0$ of $L_0$ cannot be a point since $(n-1)F\cdot L_0=2(n-2)\neq 0$. 
Therefore $\mathscr C_0$ is a curve and we obtain (i).

\begin{figure}
\includegraphics{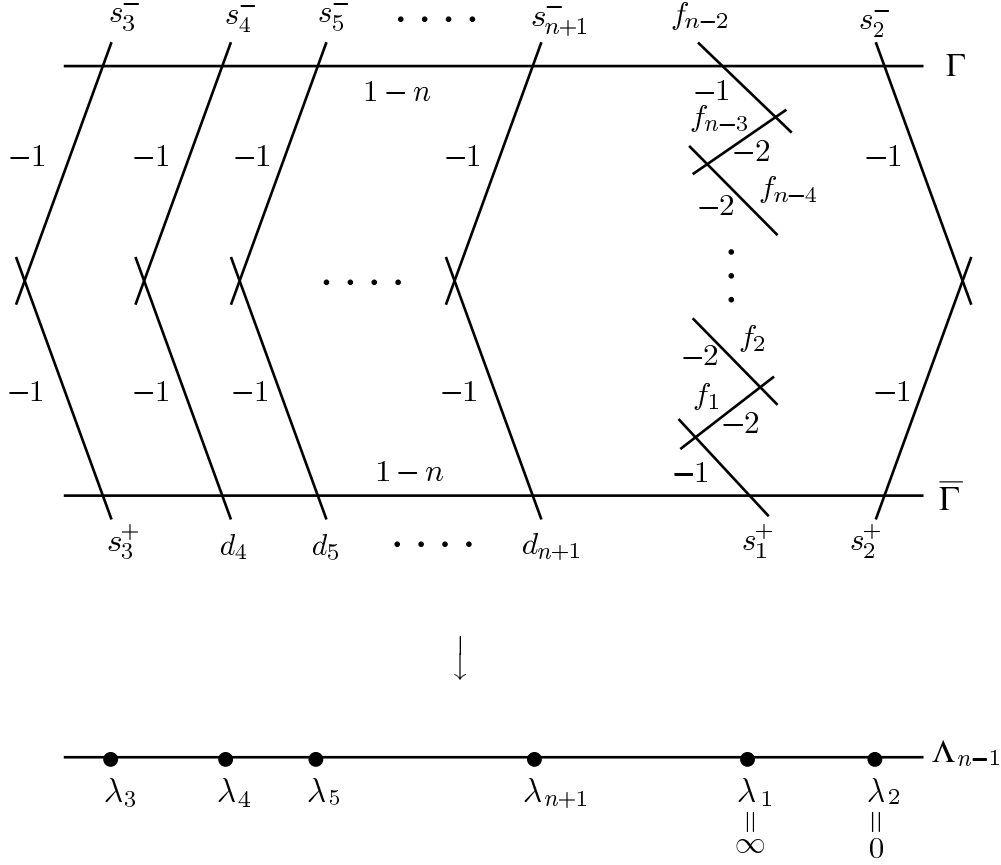}
\caption{Reducible fibers of the projection $\tilde{\mathscr T}\to\Lambda_{n-1}$}
\label{fig-mt}
\end{figure}

To show (ii) we set $\tilde{\mathscr C}_0:=\nu^{-1}(\mathscr C_0)$ and we determine the cohomology class of $\tilde{\mathscr C}_0$ on $\tilde{\mathscr T}$.
As generators of the cohomology group $H^2(\tilde{\mathscr T},\mathbf Z)\simeq\mathbf Z^{2n}$, 
we choose the following $2n$ curves:
\begin{align}\label{gen1}
\Gamma,\,f,\,s_3^+, d_i\,(4\le i\le n+1),\,s_2^-,\,f_j\,(1\le j\le n-2),
\end{align}
where  $f$ is a fiber of the conic bundle $\tilde{\mathscr T}\to\Lambda_{n-1}$, $s_i^{\pm}$ are the images of the degree-one divisors $S_i^{\pm}\subset \tilde Z$ under $\tilde{\Phi}^{\mathbf C^*}$ or $\tilde{\Phi}_{n-1}^{\mathbf C^*}$, $f_j$ and $d_i$ are the images of the exceptional divisors $F_j$ and $D_i$ respectively (arose  in obtaining $Z_n$) under the same map, and $\Gamma$ is one of the 2 exceptional curves of the partial resolution $\hat{\mathscr T}\to\mathscr T$ specified by the property that it intersects $s_2^-$. (See Figure \ref{fig-mt}.)
Then since $L_0$ intersects $S$ transversally at 2 points for general $S\in |F|$ and 
since $\tilde{\Phi}_{n-1}^{\mathbf C^*}|_S$ is holomorphic quotient map by Prop.\,\ref{prop-quot1}, we have $\tilde{\mathscr C}_0\cdot f=2$.
Similarly, we have $\tilde{\mathscr C}_0\cdot d_i=\tilde{\mathscr C}_0\cdot s_3^+=\tilde{\mathscr C}_0\cdot s_2^-=\tilde{\mathscr C}_0\cdot f_{n-2}=1$ and $\tilde{\mathscr C}_0\cdot\Gamma=\tilde{\mathscr C}_0\cdot f_j=0$ for $1\le j\le n-3$.
From this we can deduce that, as cohomology classes,
\begin{align}\label{cohom1}
\tilde{\mathscr C}_0\sim 2\Gamma+(n-1)f+s_2^-+\sum_{k=1}^{n-2}kf_k-s_3^+-\sum_{j=4}^{n+1}d_j.
\end{align}
Since we have an explicit realization of $\mathscr T$  as an embedded surface in $\mathbf{CP}^{n+1}$ as in Theorem\,\ref{thm-mt1}, we can check that the cohomology class \eqref{cohom1} is precisely the pullback of the hyperplane section class, by the minimal resolution $\nu$.
Hence we obtain (ii).
Finally, by  \eqref{cohom1} we obtain
\begin{align}
\tilde{\mathscr C}_0^2=2n-2,\,\,\,K\cdot\tilde{\mathscr C}_0=-4
\end{align}
on $\tilde{\mathscr T}$.
From this it follows that the virtual genus of $\tilde{\mathscr C}_0$ is $n-2$.
Since $\tilde{\mathscr C}_0$ and $\mathscr C_0$ are biholomorphic by $\nu$, 
the virtual genus are same.
Thus we obtain (iii).
\proofend

\begin{prop}\label{prop-discrim1}
The discriminant locus of the conic bundle $\mu(\tilde Z)\to\tilde{\mathscr T}$ consists of the following curves.
(a) The 2 exceptional curves $\Gamma$ and $\ol{\Gamma}$ of the partial resolution $\hat{\mathscr T}\to\mathscr T$.
(b) Reducible fibers of the projection $\tilde{\mathscr T}\to\Lambda_{n-1}$ over the $(n-1)$ points $(0,0\cdots,0,1)$ and $(1,\lambda_i,\lambda_i^2,\cdots,\lambda_i^{n-1})$ for $4\le i\le n+1$
in the homogeneous coordinate $(z_1,\cdots,z_n)$ on $\mathbf{CP}^{n-1}$ as before.
(c) The curve $\tilde{\mathscr C}_0=\nu^{-1}(\mathscr C_0)$, where $\mathscr C_0=\Phi_{n-1}^{\mathbf C^*}(L_0)$ as in Lemma \ref{lemma-hps1}, and $\nu:\tilde{\mathscr T}\to\mathscr T$ is the minimal resolution as before.
\end{prop}

Recall that in (b) the fiber over the point $(0,0\cdots,0,1)\in\Lambda_{n-1}$ is a string of $(n-1)$ rational curves. Other reducible fibers consist of 2 irreducible components.
Note that the proposition claims that reducible fibers of $\tilde{\mathscr T}\to\Lambda_{n-1}$ over the 2 points $(1,\lambda_i,\lambda_i^2,\cdots,\lambda_i^{n-1})$ for $i=2,3$ do not contained in the discriminant locus of $\mu(\tilde Z)\to\tilde{\mathscr T}$.

\noindent
Proof of Prop.\,\ref{prop-discrim1}. By our explicit way of the elimination for the indeterminacy locus in  neighborhoods of $C_2$ and $\ol{C}_2$, we have
$\{\Gamma,\ol{\Gamma}\}=\{\tilde{\Phi}_{n-1}^{\mathbf C^*}(E_1),\tilde{\Phi}_{n-1}^{\mathbf C^*}(\ol{E}_1)\}$.
We distinguish $\Gamma$ and $\ol\Gamma$ by supposing $\tilde{\Phi}_{n-1}^{\mathbf C^*}(E_1)=\Gamma$.
(This is compatible with the previous distinction we made in \eqref{gen1}).
Then we have $(\tilde{\Phi}_{n-1}^{\mathbf C^*})^{-1}(\Gamma)=E_1+\ol{E}_3$ (in $Z_n$), and both $E_1$ and $\ol{E}_3$ are mapped surjectively to the curve $\Gamma$ by $\tilde{\Phi}_{n-1}^{\mathbf C^*}$.
Hence $\Gamma$ is contained in the discriminant locus of $\tilde{\Phi}_{n-1}^{\mathbf C^*}$.
Since the full elimination $\tilde{\Phi}^{\mathbf C^*}:\tilde{Z}\to\tilde{\mathscr T}$ factors through $\tilde{\Phi}_{n-1}^{\mathbf C^*}$ by our choice, the inverse image $(\tilde{\Phi}^{\mathbf C^*})^{-1}(\Gamma)$ contains $E_1$ and $\ol{E}_3\,(\subset \tilde Z)$ at least.
Hence  $\Gamma$ is contained in the  discriminant locus of $\tilde{\Phi}^{\mathbf C^*}$ as well.
Further, since the blow-ups $\tilde{Z}\to Z_n$ do not touch the divisor $E_2$ by our choice, 
$E_1$ and $\ol{E}_3$ still intersect $E_2$ in $\tilde{Z}$, along  curves which are mapped biholomorphically to $\Gamma$ by $\tilde{\Phi}^{\mathbf C^*}$.
Therefore, $E_1$ and $\ol{E}_3$ are not contracted to  curves by the map $\mu$.
This means that the conic bundle $\mu(\tilde{Z})\to\tilde{\mathscr T}$ contains $\Gamma$ as a discriminant curve.
By reality, the same is true for the conjugate curve $\ol{\Gamma}$.

Next, to see that the fibers in (b) are also contained in the discriminant locus, 
we again recall that the fibers of the projection $\tilde{\mathscr T}\to\Lambda_{n-1}$ are the images of the members of the pencil $|F|$.
 In particular, the reducible  fiber over the point $(1,\lambda_i,\lambda_i^2,\cdots,\lambda_i^{n-1})\in\Lambda_{n-1}$ is the image of the member $S_i^++S_i^-$ for $2\le i\le n+1$, and the curves
 $\tilde{\Phi}_{n-1}^{\mathbf C^*}(S_i^+)$ and $\tilde{\Phi}_{n-1}^{\mathbf C^*}(S_i^-)$ are precisely the irreducible components of the fiber.
 Then we see again by the explicit way of blowing-ups that we have, for $4\le i\le n+1$,
 \begin{align}
(\tilde{\Phi}_{n-1}^{\mathbf C^*})^{-1}(\tilde{\Phi}_{n-1}^{\mathbf C^*}(S_i^+))=S_i^++D_i
\,\,\,{\text{and }}\,\,\,\,
(\tilde{\Phi}_{n-1}^{\mathbf C^*})^{-1}(\tilde{\Phi}_{n-1}^{\mathbf C^*}(S_i^-))=S_i^-+\ol{D}_i.
\end{align}
(Note that for $i=2, 3$ these do not hold because $D_2, D_3,\ol{D}_2$ and $\ol{D}_3$ do not exist.
Instead, $(\tilde{\Phi}_{n-1}^{\mathbf C^*})^{-1}(\tilde{\Phi}_{n-1}^{\mathbf C^*}(S_i^+))=S_i^+$ and $(\tilde{\Phi}_{n-1}^{\mathbf C^*})^{-1}(\tilde{\Phi}_{n-1}^{\mathbf C^*}(S_i^-))=S_i^-$ hold for $i=2,3$.
This is why the fibers over  the points $(1,\lambda_i,\lambda_i^2,\cdots,\lambda_i^{n-1})$ are not contained in the discriminant locus of $\mu(\tilde Z)\to\tilde{\mathscr T}$ for $i=2,3$.)
Hence by the same reasoning for $\Gamma$,
the fibers over $(1,\lambda_i,\lambda_i^2,\cdots,\lambda_i^{n-1})\in\Lambda_{n-1}$, $4\le i\le n+1$, are discriminant curves for $\mu(\tilde{Z})\to\tilde{\mathscr T}$.
On the other hand, the reducible  fiber over the point $(0,0\cdots,0,1)\in\Lambda_{n-1}$ corresponds to the reducible member $S_1^++S_1^-\in|F|$.
The curve $\tilde{\Phi}_{n-1}^{\mathbf C^*}(S_1^+)$ (resp.\,$\tilde{\Phi}_{n-1}^{\mathbf C^*}(F_{n-2})$) is the component of the reducible fiber which intersects $\ol{\Gamma}$ (resp.\,$\Gamma$).
We have
 \begin{align}
(\tilde{\Phi}_{n-1}^{\mathbf C^*})^{-1}(\tilde{\Phi}_{n-1}^{\mathbf C^*}(S_1^+))=S_1^++\ol{F}_{n-2}
\,\,\,{\text{and }}\,\,\,\,
(\tilde{\Phi}_{n-1}^{\mathbf C^*})^{-1}(\tilde{\Phi}_{n-1}^{\mathbf C^*}(S_1^-))=S_1^-+F_{n-2},
\end{align}
and $S_1^+$, $\ol{F}_{n-2}$, $S_1^-$ and $F_{n-2}$ are mapped surjectively to the irreducible components.
Similarly, we have 
 \begin{align}
(\tilde{\Phi}_{n-1}^{\mathbf C^*})^{-1}(\tilde{\Phi}_{n-1}^{\mathbf C^*}(F_i))=F_i+\ol{F}_{n-2-i}
\end{align}
for $1\le i\le n-3$, and $F_i$ and $\ol{F}_i$ are mapped surjectively to the components.
These mean that all the irreducible components of the fiber over $(0,0,\cdots,0,1)\in\Lambda_{n-1}$ are discriminant curves of $\tilde{\Phi}_{n-1}^{\mathbf C^*}$.
By the same reasoning for $\Gamma$ and $\ol{\Gamma}$ above, this is still true for the conic bundle  $\mu(\tilde{Z})\to\tilde{\mathscr T}$.
Thus we have seen that all fibers in (b) are actually discriminant curves.

Next we show that $\tilde{\mathscr C}_0$ is also a discriminant curve.
Since $\tilde{\Phi}_{n-1}^{\mathbf C^*}$ can be viewed as a $\mathbf C^*$-quotient map, the inverse image $(\tilde{\Phi}_{n-1}^{\mathbf C^*})^{-1}(\tilde{\mathscr C}_0)$ is a $\mathbf C^*$-invariant divisor in $Z_n$, which clearly contains the fixed twistor line $L_0$.
Further, by the explicit form of $U(1)$-action in a neighborhood of the fixed point on $n\mathbf{CP}^2$ under $L_0$,
the $U(1)$-action on the twistor space in a neighborhood of any point of $L_0$ takes the form 
\begin{align}
(u,v,w)\mapsto (su,s^{-1}v, w),\,\,s\in U(1)
\end{align}
where $L_0=\{u=v=0\}$ in the neighborhood.
Hence any $\mathbf C^*$-invariant divisor containing $L_0$ must contain at least one of the (locally defined) divisors $\{u=0\}$ and $\{v=0\}$.
Moreover, since the map $\tilde{\Phi}_{n-1}^{\mathbf C^*}$ is continuous in a neighborhood of $L_0$ (since it is holomorphic there), these two divisors must be mapped to the same curve in $\tilde{\mathscr T}$.
Hence the divisor $(\tilde{\Phi}_{n-1}^{\mathbf C^*})^{-1}(\tilde{\mathscr C}_0)$ contains both of the two divisors $\{u=0\}$ and $\{v=0\}$.
This means that $(\tilde{\Phi}_{n-1}^{\mathbf C^*})^{-1}(\tilde{\mathscr C}_0)$ has ordinary double point along $L_0$.
(Later on this will turned out to be decomposed into 2 irreducible components in $Z_n$.)
Hence $\tilde{\mathscr C}_0$ is a discriminant curve of $\tilde{\Phi}_{n-1}^{\mathbf C^*}$.
Therefore the same is true for the conic bundle 
 $\mu(\tilde{Z})\to\tilde{\mathscr T}$.
 
 Thus we have seen that all curves in (a), (b) and (c) are actually discriminant curves.
 Finally we show that there is no discriminant curve other than these.
 As is already noted,
 the discriminant curve is a member of the system $|N_{E_2/Z_n}^{\vee}\otimes N_{\ol{E}_2/Z_n}^{\vee}|$. 
Thus it is enough to show that the sum of all curves in (a), (b) and (c) already belongs to this system.
But since the normal bundles and all curves are explicitly given, this can again be verified by routine computations of intersection numbers.
\proofend

\subsection{Existence of non-trivial members and their Chern classes}\label{ss-nt}
Using Prop.\,\ref{prop-discrim1}, we can readily obtain $\mathbf C^*$-invariant divisors in the twistor spaces which will play an essential role in our analysis of the structure of the twistor spaces as follows.
(In the course of the proof, we obtain a projective model of the twistor space as a conic bundle.)

\begin{prop}\label{prop-Y1}
Let $\mathscr C_0$ be as in Lemma\,\ref{lemma-hps1}.
Then the inverse image $(\Phi_{n-1}^{\mathbf C^*})^{-1}(\mathscr C_0)$ splits into 2 irreducible components $Y$ and $\ol{Y}$  intersecting transversally along the fixed twistor line $L_0$.
Further, the degrees of $Y$ and $\ol{Y}$  are $n-1$.
\end{prop}

\noindent Proof.
First we recall that the 2 divisors $E_2$ and $\ol{E}_2$ in $\tilde Z$ are contained in the $\mathbf C^*$-fixed locus.
Since the $\mathbf C^*$-action on $\tilde Z$ is effective and non-trivial, we can suppose that  it acts on the normal bundle $N_{E_2/\tilde{Z}}$ by usual scalar multiplication on each fiber.
Hence by using the reality, the natural $\mathbf C^*$-action on the bundle 
$N_{E_2/\tilde{Z}}^{\vee}\oplus N_{\ol{E}_2/\tilde{Z}}^{\vee}\oplus\mathscr O$ is given by
\begin{align}
(x,y,t)\longmapsto (sx,s^{-1}y, t),\,\,\,s\in\mathbf C^*.
\end{align}
where $(x,y,t)$ represents points of the bundle.
Since the conic bundle $\mu(\tilde Z)$ is $\mathbf C^*$-invariant,
its defining equation in $\mathbf P(N_{E_2/\tilde{Z}}^{\vee}\oplus N_{\ol{E}_2/\tilde{Z}}^{\vee}\oplus\mathscr O)$ must be of the form
\begin{align}\label{cbdle2}
xy=P_0P_1t^2,
\end{align}
where $P_0$ is a section of a line bundle over $\tilde{\mathscr T}$ whose zero divisor is $\tilde{\mathscr C}_0$, and $P_1$ is a section of a line bundle whose zero divisor is the sum of all discriminant curves in (a) and (b) of Prop.\,\ref{prop-discrim1}.
The equation \eqref{cbdle2} immediately implies that the inverse image of the curve $\tilde{\mathscr C}_0$ splits into 2 irreducible components
$\{x=P_0=0\}$ and $\{y=P_0=0\}$, which are clearly mutually conjugate.
Hence the corresponding divisor $(\Phi_{n-1}^{\mathbf C^*})^{-1}(\mathscr C_0)$ (in $Z$) splits into two irreducible components, for which we  denote by  $Y$ and $\ol{Y}$.
Then these $Y$ and $\ol{Y}$ (in $Z$)  intersect transversally along $L_0$ since it is true already in $\mu(\tilde Z)$ and since the blow-ups $\tilde Z\to Z_n\to\cdots \to Z_1\to Z$ do not touch $L_0$.
Finally, since $\mathscr C_0$ is a hyperplane section as in Lemma \ref{lemma-hps1}, $Y+\ol{Y}$  belongs to the system $|(n-1)F|^{\mathbf C^*}$.
Hence its degree is $2(n-1)$. Hence by reality the degrees of $Y$ and $\ol{Y}$ are both $(n-1)$.
\proofend

\begin{rmk}{\em
As in the proof, the equation \eqref{cbdle2} gives a projective model of our twistor spaces as a conic bundle.
This is an analogue of the equation of the conic bundles of the twistor spaces, obtained in \cite{Hon07-3}.
A remarkable difference between the two cases is the following.
In \cite{Hon07-3}, the structure of the surface (minitwistor space) $\mathscr T$ or its resolution $\tilde{\mathscr T}$ does not depend on $n$.
Instead, the number of the irreducible components of the discriminant curve increases as $n$ does.
In contrast, for the present twistor spaces, 
the structure of $\mathscr T$ and $\tilde{\mathscr T}$ changes depending on $n$. In fact, the minimal resolution $\tilde{\mathscr T}$ has a conic bundle structure over the rational curve $\Lambda_{n-1}$, which has precisely  $(n+1)$ reducible fibers.
(As is clear from the construction, this reflects the fact that the number of reducible members of the pencil $|F|$ on the twistor spaces is $(n+1)$, which in turn is based on the number of irreducible components of the anticanonical cycle $C$. In \cite{Hon07-3}, the cycle $C$ always consists of 8 components for any $n$.)
}
\end{rmk}

Next we would like to determine the Chern classes of the divisors $Y$ and $\ol{Y}$. For this purpose we prove the following.
(Note that we have not made distinction between $Y$ and $\ol{Y}$ yet.)

\begin{lemma}\label{lemma-Y2}
Let $Y$ and $\ol{Y}$ be as in Prop.\,\ref{prop-Y1} and $S$ a real  irreducible member of the pencil $|F|$. Then the restriction $Y|_S$ or $\ol{Y}|_S$ coincides with the following curve on $S$:
\begin{align}\label{rest1}
(n-2)C_1+\sum_{i=2}^n(n+1-i)C_i\,\,+B_1+B_2.
\end{align}
\end{lemma}
\noindent
(For $B_1$ and $B_2$, see the initial construction of our surface $S$ given in \S \ref{ss-S}.) 

\noindent
Proof of Lemma \ref{lemma-Y2}.
Since $Y+\ol{Y}\in |(n-1)F|^{\mathbf C^*}$ as in the proof of Prop.\,\ref{prop-Y1} and 
Bs\,$|(n-1)F|^{\mathbf C^*}=C-C_{n+1}-\ol{C}_{n+1}$ by Prop.\,\ref{prop-generator1}, 
$Y+\ol{Y}$ contains $C_i$ and $\ol{C}_i$ for $1\le i\le n$.
We first show that the component $C_2$ (and $\ol{C}_2$ also) cannot be contained in both of $Y$ and $\ol{Y}$ simultaneously.
For this, we recall that as in the proof of Prop.\,\ref{prop-Y1}, $Y$ and $\ol Y$ are bimeromorphic images of the divisors $\{x=P_0=0\}$ and $\{y=P_0=0\}$ in the conic bundle $\mu(\tilde Z)$.
In  $\mu(\tilde Z)$, the 2 sections $E_2$ and $\ol E_2$ (of $\mu(\tilde Z)\to\tilde{\mathscr T}$) are defined by $\{x=t=0\}$ and $\{y=t=0\}$, where we have not yet specified whether $E_2=\{x=t=0\}$ or $E_2=\{y=t=0\}$ holds.
So here we suppose that  $E_2=\{x=t=0\}$ holds.
We also suppose that $Y$ and $\ol Y$ are bimeromorphic images of $\{x=P_0=0\}$ and $\{y=P_0=0\}$ respectively.
Then $\{x=P_0=0\}\cap E_2$ in $\mu(\tilde Z)$ is a curve which is biholomorphic to $\tilde{\mathscr C}_0$ by the projection $\mu (\tilde Z)\to\tilde{\mathscr T}$, and evidently $\{y=P_0=0\}\cap {E}_2=\emptyset$ in  $\mu(\tilde Z)$.
Similarly $\{x=P_0=0\}\cap \ol{E}_2=\emptyset$ and $\{y=P_0=0\}\cap \ol E_2\simeq \tilde{\mathscr C}_0$ in  $\mu(\tilde Z)$.
Then since the two bimeromorphic maps $\mu:\tilde Z\to \mu(\tilde Z)$ and $\tilde Z\to Z_n$ are biholomorphic in a neighborhood of $E_2$ and $\ol E_2$, in $Z_n$ also, $Y\cap E_2\simeq\tilde{\mathscr C}_0$ and $\ol Y\cap E_2=\emptyset$ hold.
These mean that in $Z$ also, $Y\supset C_2$ and $\ol Y\cap C_2=\emptyset$ from the above equations.

Next we show by using the computations in the proof of Lemma \ref{lemma-hps1} that the restriction $Y|_S$ contains the curve $C_2$ with multiplicity  $(n-1)$. 
For this, we recall that the exceptional divisor $E_2\subset \tilde{Z}$ is biholomorphic to the surface $\tilde{\mathscr T}$ by the map $\tilde{\Phi}^{\mathbf C^*}:\tilde Z\to\tilde{\mathscr T}$.
Let $h\in H^2(E_2,\mathbf Z)\simeq H^2(\tilde{\mathscr T},\mathbf Z)$ be the cohomology class of fibers for the natural projection $E_2\to C_2$.
Then since $Y+\ol Y=(\Phi^{\mathbf C^*}_{n-1})^{-1}(\mathscr C_0)$ and $C_2\not\subset \ol Y$,  in order to see that $Y|_S$ contains $(n-1)C_2$, 
it suffices to show that
\begin{align}\label{kemul}
h\cdot \tilde{\mathscr C}_0=n-1.
\end{align}
Using the $2n$ curves \eqref{gen1} as a basis of $ H^2(\tilde{\mathscr T},\mathbf Z)$, it can be verified that we have
\begin{align}\label{cohom2}
h\sim \Gamma+s_2^-+\sum_{k=1}^{n-2}kf_k.
\end{align}
On the other hand, the cohomology class of $\tilde{\mathscr C}_0$ is given by \eqref{cohom1}.
By \eqref{cohom1} and \eqref{cohom2} we obtain \eqref{kemul}.
Hence 
the restriction $Y|_S$ contains the curve $C_2$ with multiplicity $(n-1)$.

We next show that the divisor $(Y+\ol{Y})|_S$ contains the curve
\begin{align}\label{rest3}
(n-2)(C_1+\ol{C}_1)+\sum_{i=2}^n(n+1-i)(C_i+\ol{C}_i),
\end{align}
by using intersection numbers as a main tool.
(The curve \eqref{rest3} is exactly (\eqref{rest1}$-B_1-B_2)$ plus its conjugate curve).
Since we have
$
Y+\ol{Y}\supset C-C_{n+1}-\ol{C}_{n+1}
$ and $(Y+\ol Y)|_S\supset (n-1)(C_1+\ol C_1)$ as before, 
 the curve 
\begin{align}\label{cohom3}
\left(Y+\ol{Y}\right)|_S-\left\{
(C_1+\ol{C}_1)+(n-1)(C_2+\ol{C}_2)+\sum_{i=3}^n(C_i+\ol{C}_i)
\right\}
\end{align}
must be a zero divisor or an effective curve (on $S$).
Further, by using $(Y+\ol{Y})|_S=(n-1)K_S^{-1}$, the intersection number of the curve \eqref{cohom3} and $C_1$ can be computed to be
\begin{align}
(n-1)(2-n) +(n-1)=(n-1)(3-n).
\end{align}
Since this is negative, the curve \eqref{cohom3}  contains $C_1$.
By the same computations, we see that the curve \eqref{cohom3}  contains $C_3$ too.
Hence by reality, it  contains $\ol{C}_1+\ol{C}_3$ also.
So subtracting $(C_1+\ol{C}_1)+(C_3+\ol{C}_3)$ from \eqref{cohom3}, we obtain that 
\begin{align}\label{cohom4}
\left(Y+\ol{Y}\right)|_S-\left\{
2(C_1+\ol{C}_1)+(n-1)(C_2+\ol{C}_2)+2(C_3+\ol{C}_3)+\sum_{i=4}^n(C_i+\ol{C}_i)
\right\}
\end{align}
is still effective, or a zero divisor.
If $n=4$, this already proves the claim that $(Y+\ol{Y})|_S$ contains the curve \eqref{rest3}.
If $n\ge 5$,  the intersection numbers of the class \eqref{cohom4} with $C_1$ and $C_3$ can be computed to be $(n-1)(4-n)$ and $(4-n)$ respectively, both of which are negative.
Therefore, subtracting $(C_1+\ol{C}_1)+(C_3+\ol{C}_3)$ from \eqref{cohom4}
and then computing its intersection number with $C_4$ which turns out to be $(-2)<0$, it follows that 
the class 
\begin{align}\label{cohom5}
\left(Y+\ol{Y}\right)|_S-\left\{
3(C_1+\ol{C}_1)+(n-1)(C_2+\ol{C}_2)+3(C_3+\ol{C}_3)+2(C_4+\ol{C}_4)+\sum_{i=5}^n(C_i+\ol{C}_i)
\right\}
\end{align}
is still effective, or a zero divisor.
If $n=5$, this proves the claim that $(Y+\ol{Y})|_S$ contains the curve \eqref{rest3}.
For general $n$, by repeating this argument we obtain that 
$(Y+\ol{Y})|_S$ contains the curve \eqref{rest3}.

Next we show that the remaining curves $B_1,B_2,\ol B_1$ and $\ol{B}_2$ are also contained in $(Y+\ol{Y})|_S$.
But this is obvious since the restriction of $\Phi_{n-1}^{\mathbf C^*}:Z\to\mathscr T$ onto $S\in |F|$ is a holomorphic quotient map by Prop.\,\ref{prop-quot1}, and since both $Y$ and $\ol{Y}$ are the unions of the 2 reducible  fibers of these quotient maps consisting of  2 irreducible components.
Thus we have proved that $(Y+\ol Y)|_S$ contains the curve \eqref{rest1} and its conjugate curve.

Next we see that this inclusion is moreover an equality.
For this, by computing intersection numbers, we can verify that \eqref{rest1} plus its conjugate curve is a member of the system $|(n-1)K_S^{-1}|$.
On the other hand, we have $Y+\ol Y\in |(n-1)F|$ and $(n-1)F|S\simeq (n-1)K_S^{-1}$.
Further the restriction map $H^2(Z,\mathbf Z)\to H^2(S,\mathbf Z)$ is injective.
These imply the coincidence.

Finally we show that the curve $Y|_S$ is exactly \eqref{rest1}.
For this it is enough to see that $\ol Y$ does not contain the components $C_i$ for $4\le i\le n$.
(Recall that we have already seen that $C_1\not\subset\ol Y$ and $C_3\not\subset Y$.)
By our explicit construction of the surface $S$, we see that $\mathbf C^*$ acts on the chain $C_3+C_4+\cdots+C_{n+1}$ in such a way that if $z\in C_i$ is not a fixed point then $sz\in C_i$ goes to the fixed point $C_i\cap C_{i+1}$ as $s\to\infty$ (since $C_2$ and $\ol C_2$ are `source' and `sink'), where $C_{n+2}:=\ol C_1$.
Since $\ol Y$ actually contains $\ol C_1$, this means that if $\ol Y$ contains $C_i$ ($3\le i\le n$), then $\ol Y$ contains $C_j$ for $i\le j\le n+1$.
However, we already know that $Y+\ol Y$, and hence $\ol Y$ does not contain $C_{n+1}$.
\proofend

\bigskip
We are now ready to prove the key result, which means the existence of  `non-trivial' members of the system $|(n-2)F|$ (in the sense of Def.\,\ref{def-trivial}):

\begin{prop}\label{prop-nt4}
Let $Y$ and $\ol{Y}$ be the $\mathbf C^*$-invariant divisors as in Prop.\,\ref{prop-Y1}.
Then both of the two $\mathbf C^*$-invariant divisors
\begin{align}\label{nt11}
Y+(n-3)S_1^-\, \,\,{\text{and}}\,\,\,\,\,
\ol{Y}+(n-3)S_1^+
\end{align}
belong to the linear system $ |(n-2)F|$.
Further, if  $x_1$ and $x_2\in H^0((n-2)F)$ are mutually conjugate  sections defining divisors 
\eqref{nt11} respectively, $\mathbf C^*$  acts on these by $(x_1,x_2)\mapsto(sx_1,s^{-1}x_2)$ or  $(x_1,x_2)\mapsto(s^{-1}x_1,sx_2)$ for $s\in\mathbf C^*$.
\end{prop}

\noindent
Proof.
Since $S_1^-|_S=\sum_{i=1}^{n+1}\ol{C}_i$, 
 Lemma \,\ref{lemma-Y2} means 
\begin{align}\label{rest2}
(Y+(n-3)S_1^-)|_S=(n-2)C_1+\sum_{i=2}^n(n+1-i)C_i\,\,+B_1+B_2+(n-3)\sum_{i=1}^{n+1}\ol{C}_i.
\end{align}
It is a routine computation to verify that the right-hand side of \eqref{rest2} belongs to the system $|(n-2)K_S^{-1}|$.
(For example, it is enough to show the coincidence of their intersection numbers with the curves $\{C_i, \ol{C}_i,  B_j, \ol{B}_j\set 1\le i\le n+1,j=1,2\}$, which generate $H^2(S,\mathbf Z)$.)
Then since the restriction map $H^2(Z,\mathbf Z)\to H^2(S,\mathbf Z)$ is injective and $(n-2)F|_S\simeq (n-2)K_S^{-1}$, we obtain $Y+(n-3)S_1^-\in |(n-2)F|$.
Since $(n-2)F$ is a real bundle we also have $\ol Y+(n-3)S_1^+\in |(n-2)F|$.
Furthermore, it can be verified that the curve \eqref{rest2} is the inverse image of a $\mathbf C^*$-invariant but non-real line in $\mathbf{CP}^2$, by the morphism associated to the system $|(n-2)K_S^{-1}|$.
Therefore by Prop.\,\ref{prop-multan1} (v), precisely one of $x_1$ and $x_2$ is acted by a scalar multiplication of $\mathbf C^*$.
With the aid of the reality, this implies the final claim of the proposition.
\proofend

\section{Bimeromorphic images of the twistor spaces}
In Section 2 we considered the system $|(n-1)F|^{\mathbf C^*}$ and showed that the image of its associated map is a complex surface $\mathscr T$ whose defining equation can be explicitly determined.
The map could also be regarded as a quotient map  by $\mathbf C^*$-action
and $\mathscr T$ can be regarded as an orbit space (i.\,e.\,minitwistor space).
In Section 3 we studied this map in detail and  finally found  non-trivial members of  the system $|(n-2)F|$ (Prop.\,\ref{prop-nt4}).
Once we obtain these non-trivial members, it is possible to show that 
the map $\Phi_{n-2}$ (associated to $|(n-2)F|$)  is generically 2 to 1 covering onto its image.
However, it seems difficult yet to derive a detailed form of a defining equation of the discriminant locus of the covering.

To remedy this, in this section, we investigate the complete system $|(n-1)F|$ and show that its associate map $\Phi_{n-1}$ is bimeromorphic onto its image.
We further  give  defining equations of the image  in a projective space.
It provides another projective model of our twistor spaces which is different from the conic bundle description in Section 3.
The equations will be used in order to derive a defining equation of the branch divisor of the (generically) double covering map $\Phi_{n-2}$.


First we give  generators of the system $|(n-2)F|$  explicitly as follows.
\begin{prop}\label{prop-dc1}
Let $Z$ be a twistor space on $n\mathbf{CP}^2$ which has the complex surface $S$ constructed in \S \ref{ss-S} as a real member of the system $|F|$ as before.
Then we have the following.
(i) $\dim |(n-2)F|=n$,
(ii) As generators of the system, we can choose the following divisors.
(a) generators of the system $|V_{n-2}|$,
(b) 
$(x_1)=Y+(n-3)S_1^-$ and $(x_2)=\ol{Y}+(n-3)S_1^+$.
(iii) {\rm Bs}\,$|(n-2)F|=C-C_{n+1}-\ol{C}_{n+1}$.

\end{prop}

\noindent
Proof.
By the cohomology sequence of the sequence \eqref{ses3}, we obtain an exact sequence
\begin{align}
0\,\lra\,H^0((n-3)F)\,\lra H^0((n-2)F)\,\lra H^0((n-2)K_S^{-1}).
\end{align}
By Prop.\,\ref{prop-fundamental} and Prop.\,\ref{prop-multan1} (i) we have $H^0((n-3)F)=V_{n-3}=\mathbf C^{n-2}$.
Further by Prop.\,\ref{prop-multan1} (ii) we have $\dim H^0((n-2)K_S^{-1})=3$.
These imply $\dim H^0((n-2)F)\le n+1$.

Since the divisor $Y+(n-3)S_1^-$ is a non-trivial member,  $|V_{n-2}|$ and $Y+(n-3)S_1^-$ generate $(n-1)$-dimensional subsystem of $|(n-2)F|$.
So to prove (i) and (ii) it suffices to show that the remaining divisor $\ol{Y}+(n-3)S_1^+$ is not contained in this $(n-1)$-dimensional subsystem.
But this is obvious if we note that $\mathbf C^*$ acts trivially on $V_{n-2}$ and we have $(x_1,x_2)\mapsto (sx_1,s^{-1}x_2)$ or $(s^{-1}x_1,sx_2)$ by Prop.\,\ref{prop-nt4}, so that $x_2$ cannot be a linear combination of elements of $V_{n-2}$ and $x_1$.

For (iii), since  Bs$\,|V_{n-2}|=C$, it follows Bs\,$|(n-2)F|\subset C$.
Hence we have
\begin{equation}
{\rm Bs}\,|(n-2)F|=C\cap (Y+(n-3)S_1^-)\cap (\ol Y+(n-3)S_1^+).
\end{equation}
By Lemma  \ref{lemma-Y2}, the right-hand side is seen to be exactly $C-C_{n+1}-\ol C_{n+1}$.
\proofend

\bigskip
By using this proposition we can determine the structure of the system $|(n-1)F|$ as follows.
\begin{prop}\label{prop-bim1}
Let $Z$ be a twistor space on $n\mathbf{CP}^2$ 
as in Prop.\,\ref{prop-dc1}.
Then we have the following.
(i) $\dim |(n-1)F|=n+5$,
(ii) As generators of the system, we can choose the following divisors.
(a) generators of the $(n-1)$-dimensional system $|V_{n-1}|$,
(b) the 2   non-trivial members of $|(n-1)F|^{\mathbf C^*}$ given in Lemma \ref{lemma-nt1}.
(c) the 4 divisors defined by the following sections of $(n-1)F$:
\begin{align}
z_{n+3}:=x_1y_1,\,\,\,z_{n+4}:=x_2y_1,\,\,\,
z_{n+5}:=x_1y_2,\,\,\,z_{n+6}:=x_2y_2.
\end{align}
where $x_1$ and $x_2$ are mutually conjugate sections of $(n-2)F$ determining non-trivial members as in Prop.\,\ref{prop-nt4}.
(iii) 
{\rm Bs}\,$|(n-1)F|=C-C_{n+1}-\ol C_{n+1}$.
\end{prop}

\noindent
Proof.
From the sequence \eqref{ses1} with $m=n-1$, we obtain an exact sequence
\begin{align}
0\,\lra\,H^0((n-2)F)\,\lra H^0((n-1)F)\,\lra H^0((n-1)K_S^{-1}).
\end{align}
We have $\dim H^0((n-2)F)=n+1$ by Prop.\,\ref{prop-dc1} and $\dim H^0((n-1)K_S^{-1})=5$ by Prop.\,\ref{prop-multan2}.
These imply $\dim H^0((n-1)F)\le n+6$.
To prove (i) and (ii), it suffices to show that  $V_{n-1}=\mathbf C^n$ and $\{z_{n+i}\set 1\le i\le 6\}$ are linearly independent, where
 $z_{n+1}$ and $z_{n+2}$ denote mutually conjugate sections of $(n-1)F$ defining the 2 non-trivial members \eqref{nt1} (as in Theorem\,\ref{thm-mt1}).
We first show that $V_{n-1}$ and $z_{n+i}$, $1\le i\le 4$, generate $(n+4)$-dimensional subspace of $H^0((n-1)F)$.
By Prop.\,\ref{prop-generator1} we have $\dim (V_{n-1}+\mathbf C\langle z_{n+1},z_{n+2}\rangle)=n+2$.
From the explicitness we can verify that restriction of any member of the corresponding $(n+1)$-dimensional subsystem to $S\in |F|$ contains the curve $\ol C_2$ with multiplicity $(n-1)$.
On the other hand the restriction of the divisor $(z_{n+3})$ ($=Y+(n-3)S_1^-+S_1^++S_1^-)$ 
contains $\ol C_2$  with multiplicity only $(n-2)$.
These imply $z_{n+3}\not\in V_{n-1}\oplus \mathbf C\langle z_{n+1},z_{n+2}\rangle$.
Similarly, considering multiplicity along the component $C_2$, we obtain $z_{n+4}\not\in V_{n-1}\oplus \mathbf C\langle z_{n+1},z_{n+2},z_{n+3}\rangle$.
Hence $\dim (V_{n-1}+\mathbf C\langle z_{n+i}\set 1\le i\le 4\rangle)=n+4$.

We next show that $z_{n+5}\not\in V_{n-1}\oplus\mathbf C\langle z_{n+i}\set 1\le i\le 4\rangle$ and $z_{n+6}\not\in V_{n-1}\oplus\mathbf C\langle z_{n+i}\set 1\le i\le 5\rangle$ using $\mathbf C^*$-action.
By Prop.\,\ref{prop-nt4} we have either $(x_1,x_2)\mapsto (sx_1,s^{-1}x_2)$ or 
 $(x_1,x_2)\mapsto (s^{-1}x_1,sx_2)$ for $s\in\mathbf C^*$.
We may suppose that the former holds.
Then since $y_1$ and $y_2$ belong $H^0(F)=H^0(F)^{\mathbf C^*}$ we have
$(z_{n+3},z_{n+4},z_{n+5},z_{n+6})\mapsto (sz_{n+3},s^{-1}z_{n+4},sz_{n+5},s^{-1}z_{n+6})$ for $s\in\mathbf C^*$.
From this it readily follows  that  $z_{n+5}\in V_{n-1}\oplus\mathbf C\langle z_{n+i}\set 1\le i\le 4\rangle$ means $z_{n+5}\in \mathbf Cz_{n+3}$.
 Since $(z_{n+3})\neq(z_{n+5})$, this is a contradiction and we obtain $z_{n+5}\not\in V_{n-1}\oplus\mathbf C\langle z_{n+i}\set 1\le i\le 4\rangle$.
Similarly, $z_{n+6}\in V_{n-1}\oplus\mathbf C\langle z_{n+i}\set 1\le i\le 5\rangle$ implies $z_{n+6}\in \mathbf Cz_{n+4}$ and this is also a contradiction.
Hence we obtain $z_{n+6}\not\in V_{n-1}\oplus\mathbf C\langle z_{n+i}\set 1\le i\le 5\rangle$.
Thus we obtain (i) and (ii).

Finally since we have Bs\,$|(n-1)F|^{\mathbf C^*}=C-C_{n+1}-\ol{C}_{n+1}$ by Prop.\,\ref{prop-generator1}, we have Bs\,$|(n-1)F|\subset C-C_{n+1}-\ol{C}_{n+1}$.
Further, all the 4 divisors $(z_{n+i})$, $3\le i\le 6$, contain the cycle $C$.
These mean Bs\,$|(n-1)F|= C-C_{n+1}-\ol{C}_{n+1}$ and we obtain (iii).
\proofend

\begin{thm}\label{thm-bim}
Let $Z$ be a twistor space on $n\mathbf{CP}^2$ as in Prop.\,\ref{prop-dc1} and $\Phi_{n-1}:Z\to\mathbf{CP}^{n+5}$ the meromorphic map associated to the linear system $|(n-1)F|$.
Then $\Phi_{n-1}$ is bimeromorphic onto its image $X:=\Phi_{n-1}(Z)$.
Further, for a homogeneous coordinate $(z_1,\cdots,z_{n+6})$,  defining equations of the image $X$ are given by the following.
(a) equations in the defining ideal of the rational normal curve $\Lambda_{n-1}\subset\mathbf P^{\vee}V_{n-1}=\mathbf{CP}^{n-1}$ whose homogeneous coordinate is $(z_1,\cdots,z_n)$.
(b) the following 5 quadratic equations.
\begin{align}
z_1z_{n+5}&=z_2z_{n+3},\label{nt12}\\
z_1z_{n+6}&=z_2z_{n+4},\label{nt13}\\
z_{n+1}z_{n+2}&=z_1\{z_n-\sigma_1z_{n-1}+\sigma_2z_{n-2}-\cdots+(-1)^{n-1}\sigma_{n-1}z_1\}
,\label{nt14}\\
z_{n+3}z_{n+4}&=z_1g(z_1,\cdots,z_{n+2})\label{nt15}
\end{align}
where  $g$ is a linear polynomial of $z_1,\cdots,z_{n+2}$, and $\sigma_i$ are the elementary symmetric polynomials of $\lambda_3,\cdots,\lambda_{n+1}$ as in Theorem\,\ref{thm-mt1}.
\end{thm}

\noindent
Proof.
By the diagram \eqref{cd1} with $m=n-1$, the image $\Phi_{n-1}(Z)$ satisfies all equations in the defining ideal of $\Lambda_{n-1}\subset\mathbf P^{\vee}V_{n-1}$.
Further, since the restriction map $H^0((n-1)F)\to H^0((n-1)K_S^{-1})$ is surjective as proved in Prop.\,\ref{prop-bim1}, 
the restriction of $\Phi_{n-1}$ onto $S\in |F|$  coincides with the meromorphic map associated to $|(n-1)K_S^{-1}|$.
Since the latter map is actually (holomorphic and) birational onto its image by Prop.\,\ref{prop-multan2}, it follows from the diagram \eqref{cd1} that $\Phi_{n-1}$ is bimeromorphic onto its image.

Let $\{y_1,y_2\}$ and $\{z_1,\cdots,z_{n+2}\}$ be the same meaning as in  Theorem\,\ref{thm-mt1} and its proof.
 The latter gives a homogeneous coordinate on $\mathbf P^{\vee}H^0((n-1)F)^{\mathbf C^*}=\mathbf{CP}^{n+1}$. 
 By Theorem\,\ref{thm-mt1} these satisfy \eqref{nt14} (which is the same as \eqref{nt3}).
Further let  $\{x_1,x_2\}$ and  $\{z_{n+3},z_{n+4},z_{n+5},z_{n+6}\}$  be as in Prop.\,\ref{prop-bim1}.
Then we have
\begin{align}
z_1z_{n+5}=y_1^{n-1}\cdot x_1y_2=y_1^{n-2}y_2\cdot x_1y_1=z_2z_{n+3}
\end{align}
and we obtain \eqref{nt12}.
Similarly we have
\begin{align}
z_1z_{n+6}=y_1^{n-1}\cdot x_2y_2=y_1^{n-2}y_2\cdot x_2y_1=z_2z_{n+4}
\end{align}
and we obtain \eqref{nt13}.
Next we have
\begin{align}\label{40}
z_{n+3}z_{n+4}=x_1x_2y_1^2.
\end{align}
Now,  since $Y+\ol{Y}\in |(n-1)F|^{\mathbf C^*}$,  there is a linear polynomial $g(z_1,\cdots,z_{n+2})$ satisfying 
\begin{align}\label{g2}
\left(g\right)=Y+\ol{Y}.
\end{align}
Therefore 
since $(x_1x_2)=Y+\ol{Y}+(n-3)(S_1^++S_1^-)$ and $(y_1)=S_1^++S_1^-$, we obtain
\begin{align}
x_1x_2=g\cdot y_1^{n-3}.
\end{align}
By multiplying $y_1^2$ on both-hand sides, we obtain
\begin{align}
x_1y_1\cdot x_2y_1=gy_1^{n-1}.
\end{align}
Hence by \eqref{40} we obtain \eqref{nt15}.
\proofend

\bigskip
The structure of the bimeromorphic image $X$ becomes clearer if we blow-up $\mathbf{CP}^{n+5}$ along the center of the projection $\pi_{n-1}:\mathbf{CP}^{n+5}\to\mathbf{CP}^{n-1}$.
The center is $\mathbf{CP}^5$ which is explicitly given by $\{z_1=\cdots=z_n=0\}$ in the coordinate of Theorem \ref{thm-bim}.
By the blowing-up, $\mathbf{CP}^{n+5}$ becomes biholomorphic to the total space of the $\mathbf{CP}^6$-bundle $\mathbf P(\mathscr O(1)^{\oplus 6}\oplus \mathscr O)\to\mathbf{CP}^{n-1}$.
Restricting over $\Lambda_{n-1}$ whose degree is $(n-1)$, we obtain the bundle $\mathbf P(\mathscr O(n-1)^{\oplus 6}\oplus \mathscr O)\to\Lambda_{n-1}=\mathbf{CP}^1$.
Let $\hat X$ be the strict transform of $X=\Phi_{n-1}(Z)$ into this $\mathbf{CP}^6$-bundle.
Putting $\xi_i=z_{n+i}/z_1$ for $1\le i\le 6$ and using $(\xi_1,\cdots,\xi_6)\in\mathscr O(1)^{\oplus 6}$ as a non-homogeneous fiber coordinate over $\Lambda_{n-1}\backslash\{(0,\cdots,0,1)\}$, the equation of $\hat X$ is given by
\begin{align}
\xi_5&=\lambda\xi_3,\label{nt16}\\
\xi_6&=\lambda\xi_4,\label{nt17}\\
\xi_1\xi_2&=\lambda (\lambda-\lambda_3)(\lambda-\lambda_4)\cdots(\lambda-\lambda_{n+1}),\label{nt18}\\
\xi_3\xi_4&=g(1,\lambda,\cdots,\lambda^{n-1},\xi_1,\xi_2),\label{nt19}
\end{align}
where $\lambda=z_2/z_1=y_2/y_1$ is a non-homogeneous coordinate on $\Lambda_{n-1}\backslash\{(0,\cdots,0,1)\}$ as in Theorem\,\ref{thm-mt1}.
The 2 equations \eqref{nt16} and \eqref{nt17} determine a $\mathbf{CP}^4$-subbundle in the $\mathbf{CP}^6$-bundle.
In each fiber of this $\mathbf{CP}^4$-bundle,  \eqref{nt18} and \eqref{nt19} determine a quartic surface.
If $\lambda\in\mathbf R$ and if $\lambda\neq\lambda_i$ for $2\le i\le n+1$ (recall $\lambda_2=0$), 
this quartic surface is nothing but the birational image of the corresponding real irreducible member $S\in|F|$, by the map associated to $|(n-1)K_S^{-1}|$.
Thus, the bimeromorphic image $X=\Phi_{n-1}(Z)$ is bimeromorphic to a fiber space over $\Lambda_{n-1}=\mathbf{CP}^1$ whose general fibers are irreducible quartic surfaces.
If $\lambda=\lambda_i$ for $2\le i\le n+1$ or $\lambda=\infty$, then the right-hand side of \eqref{nt18} vanishes and consequently the fiber degenerates into 2 irreducible components.
Of course, these are the images of reducible members $S_i^++S_i^-$ for $1\le i\le n+2$.

\section{Projective models as generically double coverings}
In this section by using the results in the previous section we investigate  the meromorphic map $\Phi_{n-2}$ (associated to the system $|(n-2)F|$ on the twistor space) and show that it gives a  generically 2 to 1 covering onto its image.
Next we derive a defining equation of the discriminant locus of the double covering by using defining equations of the image $X=\Phi_{n-1}(Z)$ obtained in Theorem \ref{thm-bim}.

In order to investigate a relation between two meromorphic maps  $\Phi_{n-1}$ and $\Phi_{n-2}$, we consider an injection $H^0((n-2)F)\ra H^0((n-1)F)$ given by $\zeta\mapsto \zeta\otimes y_1$, where $y_1\in H^0(F)$ satisfies $(y_1)=S_1^++S_1^-$ as before.
Let 
\begin{align}
f:\mathbf P^{\vee}H^0((n-1)F)=\mathbf{CP}^{n+5}\,\lra\,\mathbf P^{\vee}H^0((n-2)F)=\mathbf{CP}^{n}
\end{align}
be the projection induced from the injection.
Recall that by Prop.\,\ref{prop-dc1} we can choose a set of sections
\begin{align}
y_1^{n-2},\,y_1^{n-3}y_2,\cdots,y_2^{n-2},\, x_1,\,x_2
\end{align}
as a basis of $H^0((n-2)F)\simeq\mathbf C^{n+1}$.
By taking (tensor) products with $y_1$ we obtain a set of sections 
\begin{align}
z_1=y_1^{n-1},\,z_2=y_1^{n-2}y_2,\cdots,z_{n-1}=y_1y_2^{n-2},\, z_{n+3}=x_1y_1,\,z_{n+4}=x_2y_1
\end{align}
which is a part of a basis of $H^0((n-1)F)$ given in Prop.\,\ref{prop-bim1}.
Hence the projection $f$ is explicitly given by
\begin{align}\label{proj1}
(z_1,\cdots,z_n, z_{n+1},\cdots,z_{n+6})\longmapsto (z_1,\cdots,z_{n-1}, z_{n+3},z_{n+4}).
\end{align}
Combining various meromorphic maps appeared so far, we obtain the following commutative diagram
\begin{equation}\label{cd-big}
\xymatrix{
\mathbf{CP}^{n+5}\ar[rr]^{\pi_{n-1}}\ar[dd]_f& &\mathbf{CP}^{n-1}\ar[dd]\\
&Z\ar[lu]^{\Phi_{n-1}}\ar[ru]_{\Psi_{n-1}}\ar[ld]_{\Phi_{n-2}}\ar[rd]^{\Psi_{n-2}}&\\
\mathbf{CP}^{n}\ar[rr]_{\pi_{n-2}} & &\mathbf{CP}^{n-2}
}
\end{equation}
where $\Psi_{n-1}(Z)=\Lambda_{n-1}$ and $\Psi_{n-2}(Z)=\Lambda_{n-2}$.
Of course the right projection $\mathbf{CP}^{n-1}\to\mathbf{CP}^{n-2}$ is the one induced by the inclusion $V_{n-2}\to V_{n-1}$ given by $\zeta\mapsto\zeta\otimes y_1$,
and it is explicitly given by $(z_1,\cdots,z_n)\mapsto(z_1,\cdots,z_{n-1})$.
This maps the curve $\Lambda_{n-1}$ isomorphically to $\Lambda_{n-2}$.
By blowing-up along the centers of $\pi_{n-1}$ and $\pi_{n-2}$, and then restricting over $\Lambda_{n-1}$ and $\Lambda_{n-2}$, we obtain the diagram
\begin{equation}\label{cd7}
 \CD
\mathbf{P}(\mathscr O(n-1)^{\oplus 6}\oplus\mathscr O)@>>>\Lambda_{n-1} \\
 @V{\hat f}VV \hspace{-5mm}@VVV\\
\mathbf{P}(\mathscr O(n-2)^{\oplus 2}\oplus\mathscr O)@>>>\Lambda_{n-2}\\
 \endCD
 \end{equation}
 where $\hat f$ is the meromorphic map induced from $f$.
The indeterminacy locus of  $\hat f$ is $\mathbf{CP}^3$-subbundle of the $\mathbf{CP}^6$-bundle.
This locus actually intersects the strict transform $\hat X$ of $X$  (as in the final part of the previous section).
Hence $\hat f$ has indeterminacy locus on $\hat X$.
If we blow-up along the $\mathbf{CP}^3$-subbundle, we obtain a morphism $\tilde f:\tilde X\to\mathbf{P}(\mathscr O(n-2)^{\oplus 2}\oplus\mathscr O)$, where $\tilde X$ denotes the strict transform of $\hat X$ into the blown-up space.

 We put 
 \begin{equation}\label{eta}
 \eta_1=z_{n+3}/z_1,\,\,\,\eta_2=z_{n+4}/z_1
 \end{equation}
  and use $(\eta_1,\eta_2)$ as a non-homogeneous fiber coordinate on the $\mathbf{CP}^2$-bundle in \eqref{cd7}.
 
We have now reached the main result of this paper.

\begin{thm}\label{thm-dc2}
Let $Z$ be a twistor space on $n\mathbf{CP}^2$ as in Prop.\,\ref{prop-dc1} and $\Phi_{n-2}:Z\to\mathbf{CP}^n$ the meromorphic map associated to the system $|(n-2)F|$ as before.
Then we have the following.
 (i)  $\Phi_{n-2}(Z)=(\pi_{n-2})^{-1}(\Lambda_{n-2})$ holds. 
In other words, if $\hat X$ denotes the strict transform of the bimeromorphic image $X=\Phi_{n-1}(Z)$,
we have  $\hat f(\hat X)=\mathbf{P}(\mathscr O(n-2)^{\oplus 2}\oplus\mathscr O)$.
(ii) The surjective morphism $\tilde f:\tilde X\to\mathbf{P}(\mathscr O(n-2)^{\oplus 2}\oplus\mathscr O)$ (obtained above) is generically 2 to 1.
(iii) In the above non-homogeneous coordinate $(\eta_1,\eta_2)$ on the $\mathbf{CP}^2$-bundle, the defining equation of the branch divisor $B$ of $\tilde f$ is given by the following equation.
\begin{align}\label{branch1}
\left\{
\eta_1\eta_2-\hat{g}(\lambda)
\right\}^2
=\lambda(\lambda-\lambda_3)(\lambda-\lambda_4)\cdots(\lambda-\lambda_{n+1})
\end{align}
where $\hat{g}(\lambda)$ is a polynomial with real coefficients whose degree is at most $(n-1)$.
\end{thm} 

\begin{rmk}\label{rmk-3}
{\em
As we will see below, the polynomial $\hat{g}$ in \eqref{branch1} is obtained from the polynomial $g$ in \eqref{nt15} by putting
\begin{equation}\label{dcp5}
g(1,\lambda,\cdots,\lambda^{n-1},\xi_1,\xi_2)=\hat g(\lambda)+c\xi_1+\ol{c}\xi_2,
\end{equation}
where $c$ is a non-zero constant. 
}
\end{rmk}

\begin{rmk}\label{rmk-4}
{\em
The image $\Phi_{n-2}(Z)=(\pi_{n-2})^{-1}(\Lambda_{n-2})$ is a rational scroll of planes in $\mathbf{CP}^n$ and has a cyclic quotient singularities along the center (= a line) of the projection $\pi_{n-2}$, where the order of the cyclic group is $n-2$.
}
\end{rmk}

\begin{rmk}
{\em
If $n=3$, the equation \eqref{branch1} becomes 
\begin{align}
\left\{
\eta_1\eta_2-\hat{g}(\lambda)
\right\}^2
=\lambda(\lambda-\lambda_3)(\lambda-\lambda_4)
\end{align}
where $\hat g$ is a quadratic polynomial.
In effect, this is exactly the equation of the branch quartic surface we  obtained in \cite{Hon07-2}.
In this sense, the present twistor spaces can be regarded as a generalization of the twistor spaces on $3\mathbf{CP}^2$ of `double solid type' into $n\mathbf{CP}^2$, $n$ arbitrary.
}
\end{rmk}

\noindent Proof of Theorem \ref{thm-dc2}.
Since the restriction map $H^0((n-2)F)\to H^0((n-2)K_S^{-1})$ is surjective as proved in Prop.\,\ref{prop-dc1}, the restriction of $\Phi_{n-2}$ onto $S\in|F|$ coincides with the rational map associated to the system $|(n-2)K_S^{-1}|$.
By Prop.\,\ref{prop-multan1} the latter map gives a generically 2 to 1 covering 
onto $\mathbf{CP}^2$.
Hence by the diagram \eqref{cd1} with $m=n-2$ we obtain that the map $\hat f:\hat X\to\mathbf{P}(\mathscr O(n-2)^{\oplus 2}\oplus\mathscr O)$ is surjective and generically 2 to 1.
These mean (i) and (ii).

Next to prove (iii) we use the non-homogeneous fiber coordinates $(\xi_1,\cdots,\xi_6)$ and $(\eta_1,\eta_2)$ on 
the bundles $\mathbf P(\mathscr O(n-1)^{\oplus 6}\oplus \mathscr O)$ and 
$\mathbf P(\mathscr O(n-2)^{\oplus 2}\oplus \mathscr O)$ respectively.
In these coordinates, by \eqref{proj1}, the projection $\hat f$ is explicitly given by
\begin{align}\label{proj2}
\hat f:(\xi_1,\xi_2,\xi_3,\xi_4,\xi_5,\xi_6)\longmapsto(\eta_1,\eta_2)=(\xi_3,\xi_4).
\end{align}
We recall that  defining equations of $\hat X$ are explicitly given by \eqref{nt16}--\eqref{nt19}.
Let $(\eta_1,\eta_2)$ be a point. Then by \eqref{proj2}, a point $(\xi_1,\cdots,\xi_6)$ belongs to $\hat f^{-1}(\eta_1,\eta_2)$ iff $\xi_3=\eta_1$ and $\xi_4=\eta_2$ hold.
Moreover by \eqref{nt16}--\eqref{nt19} if the point belongs to $\hat X$ further, we have
\begin{align}
\xi_1\xi_2&=\lambda (\lambda-\lambda_3)(\lambda-\lambda_4)\cdots(\lambda-\lambda_{n+1}),\label{nt20}\\
\eta_1\eta_2&=g(1,\lambda,\cdots,\lambda^{n-1},\xi_1,\xi_2).\label{nt21}
\end{align}
Therefore the point $(\eta_1,\eta_2)$ belongs to the branch point of $\tilde f:\tilde X\to\mathbf P(\mathscr O(n-2)^{\oplus 2}\oplus \mathscr O)$ iff the equations \eqref{nt20} and \eqref{nt21} have a unique solution, viewed as equations for $(\xi_1,\xi_2)$.
If we write $g$ as in  \eqref{dcp5}, \eqref{nt21} is equivalent to $\xi_2=\ol{c}^{-1}\{\eta_1\eta_2-\hat g(\lambda)-c\xi_1\}.$ 
Substituting this into \eqref{nt20} and multiplying $\ol c$ to its both-hand sides, we obtain
\begin{align}
c\xi_1^2-\{\eta_1\eta_2-\hat g(\lambda)\}\xi_1+\ol{c}\lambda  (\lambda-\lambda_3)(\lambda-\lambda_4)\cdots(\lambda-\lambda_{n+1})=0
\end{align}
This has a unique solution iff the discriminant vanishes:
\begin{align}\label{disc5}
\{\eta_1\eta_2-\hat g(\lambda)\}^2-4|c|^2\lambda  (\lambda-\lambda_3)(\lambda-\lambda_4)\cdots(\lambda-\lambda_{n+1})=0.
\end{align}
Now recalling that $g(z_1,\cdots,z_{n+2})$ was originally a section of the line bundle $(n-1)F$ whose zero divisor is $Y+\ol{Y}$ as in \eqref{g2}, we may suppose that $|c|=1/2$ by multiplying a non-zero constant.
Thus the equation \eqref{disc5} becomes
 \begin{align}\label{disc6}
\{\eta_1\eta_2-\hat g(\lambda)\}^2=\lambda (\lambda-\lambda_3)(\lambda-\lambda_4)\cdots(\lambda-\lambda_{n+1}).
\end{align}
Hence we obtain (iii).
\proofend
 
 \bigskip  
 The structure of the branch divisor $B$ in Theorem \ref{thm-dc2} is described as follows.

\begin{prop}\label{prop-B}
Let $B$ be the branch divisor  defined by the equation \eqref{branch1}.
Then we have the following.
(i) $B$ is irreducible and birational to a ruled surface of genus $[(n-1)/2]$, where 
$[k]$ denotes the biggest integer not greater than $k$.
(ii) Fibers of the natural projection $B\to\mathbf{CP}^1$ (induced from the projection $\mathbf{P}(\mathscr O(n-2)^{\oplus 2}\oplus\mathscr O)\to\mathbf{CP}^1)$ is 
non-reduced iff $\lambda=\lambda_i$ for $2\le i\le n+1$ or $\lambda=\infty$.
Moreover, the surface $B$ has rational double points of type $A_{3n-9}$ along the fiber over $\lambda=\infty$ which is a sum two lines in $\mathbf{CP}^2$.
\end{prop}
 
 Since our equation \eqref{branch1} of $B$ is explicit, it is not difficult to derive the conclusions. We leave it to the interested reader. Note that (ii) means that if $n\ge 4$,  $B$ is non-normal.

 \bigskip
 Next we show that the polynomial $\hat g(\lambda)$ in the defining equation \eqref{branch1} of the branch divisor must satisfy certain constraint.
 It is a generalization of a constraint appeared in \cite[Prop.\,2.3 or Condition (A)]{Hon07-2}, which was expressed in terms of a double root of a polynomial.
 
 \begin{prop}\label{prop-dr}
 Let $Z$ be a twistor space on $n\mathbf{CP}^2$  as in Theorem \ref{thm-dc2}, 
 whose projective model $\hat X$ is realized as a double covering of the bundle $\mathbf P(\mathscr O(n-2)^{\oplus 2}\oplus \mathscr O)\to\mathbf{CP}^1$ branched along a surface \eqref{branch1}. 
Consider  the double covering of $\mathbf{CP}^1$ whose branch locus is given by
\begin{align}\label{branch2}
\hat g(\lambda)^2-\lambda (\lambda-\lambda_3)(\lambda-\lambda_4)\cdots(\lambda-\lambda_{n+1})=0.
\end{align}
Then this double covering is a rational curve.
 \end{prop}
 
 \noindent
 Proof.
 Let $L_0\subset Z$ be the fixed line as before.
 Then since Bs\,$|(n-2)F|=C-C_{n+1}-\ol{C}_{n+1}$ by Prop.\,\ref{prop-dc1},
the meromorphic map  $\Phi_{n-2}$ is holomorphic on $L_0$.
Further, since $(n-2)F\cdot L_0=2(n-2)\neq 0$, the image $\Phi_{n-2}(L_0)$ cannot be a point.
Hence $\Phi_{n-2}(L_0)$ is a $\mathbf C^*$-fixed curve in $\mathbf{CP}^n$.
Let $l_0$ be the image curve of $L_0$ and $\Phi_{n-2}(L_0)$ into the bundle $\mathbf P(\mathscr O(n-2)^{\oplus 2}\oplus \mathscr O)\to\mathbf{CP}^1$.
Since $\mathbf C^*$ acts on fibers of this bundle as $(\eta_1,\eta_2)\mapsto (s\eta_1,s^{-1}\eta_2)$ or $(s^{-1}\eta_1,s\eta_2)$ by the choice of $(\eta_1,\eta_2)$ in \eqref{eta},
$\mathbf C^*$-fixed locus of this bundle consists of 3 sections.
But two of them are conjugate pair.
Therefore we obtain $l_0=\{\eta_1=\eta_2=0\}$.
Namely the `zero section' of the $\mathbf{CP}^2$-bundle is exactly the image of the fixed twistor line $L_0$.
Thus we obtain a holomorphic map from $L_0$ to $l_0$, which preserves the real structure.
The real structure on $L_0$ has no real point and that on $l_0$ has real point,
since the real structure acts on the base space $\Lambda_{n-2}=\mathbf{CP}^1$ by complex conjugation.
Therefore the map $L_0\to l_0$ cannot be isomorphic.
Hence, since $\Phi_{n-2}$ is generically 2 to 1, the map $L_0\to l_0$ is 2 to 1.
The branch locus of this map is the intersection of $l_0$ and  the branch locus \eqref{branch1} of $\Phi_{n-2}$.
It is exactly \eqref{branch2}.
Since $L_0$ is of course rational, this implies the claim of the proposition.
\proofend
 
 \bigskip
Of course, Prop.\,\ref{prop-dr} means that the equation \eqref{branch2} has multiple roots.
Next by using Theorem \ref{thm-dc2} and this proposition we compute the dimension of the moduli space of our twistor spaces.
For this we first count the number of parameters contained in the equation \eqref{branch1} of the branch divisor $B$.
 Since  $\deg\hat g=n-1$ in general,  $\hat g$ contains $n$ real parameters.
 (Note that we lost freedom of multiplying non-zero constants when we set $|c|=1/2$ in \eqref{disc5}.)
Moreover, there are $n-1$ real parameters $\lambda_3,\lambda_4,\cdots,\lambda_{n+1}$ contained.
On the other hand, for the non-homogeneous coordinate $\lambda$ on $\mathbf{CP}^1$,  the coordinate change $\lambda\mapsto c\lambda$, $c\in\mathbf R^*$ is allowed.
   This drops the dimension by one.
 Further, the constraint obtained in Prop.\,\ref{prop-dr} drops the number of parameters by $n-2$.
 To see this, we note that since the degree of the left-hand side of \eqref{branch2} is $2(n-1)$ in general, the virtual genus of the double covering of $\mathbf{CP}^1$ branched at the roots of \eqref{branch2} is $n-2$.
Hence the rationality of the double cover drops the dimension by $n-2$.
In conclusion, the dimension of the moduli space of our twistor spaces becomes
\begin{equation}\label{moduli1}
\{n+(n-1)\}-\{1+(n-2)\}=n.
\end{equation}

 \bigskip
 Finally we make a remark on the existence of our twistor spaces.
 Let $Z_{\rm LB}$ be a LeBrun twistor space on $n\mathbf{CP}^2$, which admits not only $\mathbf C^*$-symmetries but also $(\mathbf C^*)^2$-symmetries.
 Namely the associated LeBrun metric is supposed to admit not only $U(1)$-action but also $U(1)^2$-action.
In the paper \cite{Hon07-1} we studied $U(1)$-equivariant small deformations of the twistor space $Z_{\rm LB}$.
In particular, we have determined which $U(1)$-subgroup of $U(1)^2$ admits equivariant deformation whose resulting twistor spaces are not LeBrun twistor spaces.
The result (Prop.\,2.1 of \cite{Hon07-1}) says that there are precisely $(n-1)$ subgroups $K_i\subset U(1)^2$, $1\le i\le n-1$, which satisfy this property.
(These subgroups are specified in terms of which irreducible component is fixed by the $K_i$-action, among the anticanonical cycle in a smooth toric surface contained in $Z_{\rm LB}$.)
Among these $(n-1)$ subgroups, $K_1$-equivariant deformations of $Z_{\rm LB}$ can yield the present twistor spaces we have studied in this paper.
To see this, since all results in this paper rely on the structure of the surface $S$ contained in the system $|F|$, it suffices to verify that a smooth toric surface $S_{\rm LB}$ contained in $|F|$ of $Z_{\rm LB}$ can be $K_1$-equivariantly deformed into our surface $S$ constructed in \S \ref{ss-S}, and that the divisor $S_{\rm LB}$ survives under $K_1$-equivariant deformations of $Z_{\rm LB}$.
These properties can be proved by the same argument we have given in \cite[\S 5.1]{Hon07-3} in proving the existence of the twistor spaces studied in the paper.
Here we only remark that the dimensions of the moduli space actually coincide:
in Prop.\,2.1 of \cite{Hon07-1} we have shown that the moduli space of  non-LeBrun self-dual metrics on $n\mathbf{CP}^2$ obtained by $K_1$-equivariant deformation (of LeBrun metric with torus action) is $n$-dimensional, by determining $U(1)^2$-action on the cohomology group $H^1(\Theta_{Z_{\rm LB}})$ governing small deformations of $Z_{\rm LB}$.
This is (of course) equal to the dimension we have obtained in \eqref{moduli1}.
(See also \cite[Example 2.4]{Hon07-1}.)

\small
\vspace{13mm}
\hspace{7.5cm}
$\begin{array}{l}
\mbox{Department of Mathematics}\\
\mbox{Graduate School of Science and Engineering}\\
\mbox{Tokyo Institute of Technology}\\
\mbox{2-12-1, O-okayama, Meguro, 152-8551, JAPAN}\\
\mbox{{\tt {honda@math.titech.ac.jp}}}
\end{array}$

\end{document}